\documentclass[11pt,reqno]{amsart}
\usepackage{amsmath,amssymb,amscd,amsthm}
\usepackage[pagewise]{lineno}

\usepackage{xcolor}
\usepackage{verbatim}
\usepackage{tikz-cd}
\usetikzlibrary{babel}
\usetikzlibrary{decorations.pathmorphing}
\usepackage[margin=1.3in]{geometry}
\usepackage{hyperref}

\hypersetup{hidelinks, colorlinks=true, allcolors=blue, pdfstartview=Fit, breaklinks=true, pdftitle={Prefix exchanges and approximate K-conjugacy for foldings of de Bruijn graphs}, pdfauthor={Sihan Wei}}
\usepackage[all]{xy}
\DeclareMathOperator{\Ad}{Ad}
\usepackage{enumerate}
\usepackage{paralist}
\usepackage{amscd}
\IfFileExists{extarrows.sty}{\usepackage{extarrows}}{}

\newtheorem{thm}{Theorem}[subsection]
\newtheorem{lem}[thm]{Lemma}
\newtheorem{prop}[thm]{Proposition}

\newtheorem{cor}[thm]{Corollary}

\theoremstyle{definition}\newtheorem{df}[thm]{Definition}
\theoremstyle{definition}\newtheorem{rem}[thm]{Remark}
\newtheorem{exm}[thm]{\it Example}

\theoremstyle{definition}

\renewcommand{\phi}{\varphi}

\newcommand{\Z}{\mathbb{Z}}

\newcommand{\C}{\mathbb{C}}
\newcommand{\T}{\mathbb{T}}

\numberwithin{equation}{section}

\newcommand{\id}{\operatorname{id}}

\newcommand{\beq}{\begin{eqnarray}}
\newcommand{\eneq}{\end{eqnarray}}

\makeatletter

\def\@tocline#1#2#3#4#5#6#7{%
  \relax
  \ifnum #1>\c@tocdepth
  \else
    \par
    \addpenalty\@secpenalty
    \addvspace{#2}%
    \begingroup
      \hyphenpenalty\@M
      \@ifempty{#4}{%
        \@tempdima\csname r@tocindent\number#1\endcsname\relax
      }{%
        \@tempdima#4\relax
      }%
      \parindent\z@
      \leftskip#3\relax
      \rightskip\@pnumwidth
      \parfillskip-\@pnumwidth
      #5%
      \leavevmode
      #6%
      \nobreak
      \leaders\hbox{\normalfont.\kern0.5em}\hfill
      \nobreak
      \hbox to\@pnumwidth{\hss\normalfont #7}%
      \par
    \endgroup
  \fi
}

\renewcommand{\l@section}{%
  \@tocline{1}{0.8em}{0pt}{2.5em}{\bfseries}%
}

\renewcommand{\l@subsection}{%
  \@tocline{2}{0pt}{2.5em}{3.5em}{}%
}

\makeatother

\newcommand{\Od}{\mathcal O_d}
\newcommand{\OG}{\mathcal O_G}
\newcommand{\fullshift}{\Sigma_d^+}

\newcommand{\proofstep}[2]{\par\smallskip\noindent
  \textbf{\upshape(\ref{#1}) \(\boldsymbol{\Rightarrow}\) (\ref{#2}).}\enspace\ignorespaces}
\newtheorem*{mainA}{Theorem A}
\newtheorem*{mainB}{Theorem B}
\begin{document}
\title[Prefix exchanges and approximate $K$-conjugacy]{Prefix exchanges and approximate $K$-conjugacy\\ for foldings of de Bruijn graphs}
\author{Sihan Wei}
\date{}
\begin{abstract}
Let $G$ be the quotient of a de Bruijn graph by a folding on $d\ge3$
letters. We prove that the one-sided edge shift of $G$ is approximately
$K$-conjugate to the full $d$-shift if and only if every vertex of $G$
has indegree congruent to $1$ modulo $d-1$. The same congruences
characterize when any two members of the clopen partition determined
by the synchronizing map can be interchanged by a transposition in the
Higman--Thompson group $V_{d,1}$. 

For two foldings on the same alphabet
with $s$ and $t$ vertices, we show that the finite coefficient
$K$-condition holds precisely when $\gcd(s,d-1)=\gcd(t,d-1)$.
We also construct a two-vertex shift of entropy $\log 2$ which is
approximately $K$-conjugate to the full $3$-shift.
\end{abstract}
\maketitle

\section{Introduction}\label{sec:intro}

For a one-sided shift of finite type $(X,\sigma_X)$, how much of the
action of $\sigma_X$ can be recovered from the associated $C^*$-algebra
and the inclusion of $C(X)$? Cuntz and Krieger associated a
$C^*$-algebra $\mathcal O_A$ to a topological Markov shift with
transition matrix $A$ \cite{ck}. The path cylinders give projections
in $\mathcal O_A$ which generate its canonical diagonal. For more
general one-sided subshifts, Carlsen's Cuntz--Pimsner construction
\cite{carlsencp} and the groupoid description of Brix and Carlsen
\cite{brixcovers} give a canonical embedding
$j_X:C(X)\to\mathcal O_X$. For one-sided irreducible shifts of finite
type which are not finite periodic orbits, Matsumoto and Matui
characterized continuous orbit equivalence by the isomorphism class of
the Cuntz--Krieger algebra and the sign of $\det(I-A)$
\cite[Theorem~3.6]{matmat}. Brix and Carlsen characterized one-sided
conjugacy by a diagonal-preserving isomorphism intertwining a canonical
completely positive map \cite[Theorem~3.3]{brixconj}. Here we ask what
dynamical information is retained when compatibility with the
embeddings $j_X$ is imposed in $K$-theory.

For minimal Cantor systems, Giordano, Putnam and Skau classified
strong orbit equivalence by the ordered $K_0$-group of the crossed
product with its distinguished order unit \cite{gps}. Lin and Matui
showed that this invariant also classifies approximate $K$-conjugacy
\cite[Definition~5.3 and Theorem~5.4]{linmatui}. For classically rigid
minimal Cantor--circle extensions, Lin, Liu and Wei formulated the
finite coefficient $K$-condition in terms of the maps induced by the
coefficient embeddings on finitely generated subgroups \cite{llw}.
For the one-sided shifts considered here, $K_0(C(X))=C(X,\mathbb Z)$.
Thus the condition can be examined on characteristic functions of
clopen sets. More precisely, for two shifts $X,Y$ it requires a
unit-preserving graded isomorphism on $K$-theory and homeomorphisms
in each direction which eventually make the diagrams for $j_X,j_Y$
commute on each prescribed finitely generated subgroup
(Definition~\ref{def:finitecoeff}). In approximate $K$-conjugacy the
homeomorphisms must approximately intertwine the shifts, and the
algebra isomorphisms, with fixed inverse $KL$-classes, must approximate
the action of those same homeomorphisms on the coefficient algebras
(Definition~\ref{def:ak}).

Consider a finite alphabet $\mathcal A$ of cardinality $d\ge3$.
In the de Bruijn graph of order $L\ge1$, the vertices are the words
of length $L$, and an edge deletes the first letter and appends a new
one. A folding is an equivalence relation on the vertices preserved
by each of these transitions. In the quotient graph $G$, reading any
$u\in\mathcal A^L$ from any initial vertex leads to the same vertex
$t(u)$. In the terminology of Bleak, Cameron and Olukoya, $t$ is the
\emph{synchronizing map} and $t(u)$ is the state forced by $u$
\cite[Section~3.3 and Proposition~3.4]{bleak}. Writing
$\fullshift=\mathcal A^{\mathbb N_0}$ and $[u]$ for the cylinder with
prefix $u$, we obtain a clopen partition
\[
 P_v=\bigsqcup_{t(u)=v}[u]\quad (v\in G^0)
\]
of the full shift. The Higman--Thompson group $V_{d,1}$ acts on
$\fullshift$ by finite prefix replacements \cite{higman,matuifull}.
For disjoint nonempty clopen sets, Matui's theorem gives a
transposition interchanging them precisely when their indicators have
the same class in the zeroth homology group of the shift groupoid
\cite[Theorem~6.4 and Section~6.4]{matuifull}. For $P_v$ and $P_w$,
this says that the numbers of length-$L$ words forcing $v$ and $w$
are congruent modulo $d-1$.

\medskip\noindent\textbf{Question 1.}
Can one decide from the quotient graph $G$ whether every pair
$P_v,P_w$ admits such a transposition? What does this property imply
for the one-sided edge shift $X_G$?

\begin{mainA}
Let $G$ be the quotient by a folding of a de Bruijn graph on $d\ge3$
letters. The following are equivalent.
\begin{enumerate}
\item\label{a:ak} $X_G$ and $\fullshift$ are approximately $K$-conjugate.
\item\label{a:ind} Each vertex of $G$ has indegree congruent to $1$ modulo $d-1$.
\item\label{a:ex} For each pair $v\ne w$, a transposition in $V_{d,1}$
interchanges $P_v$ and $P_w$ and fixes the complement of their union
pointwise.
\end{enumerate}
\end{mainA}

To pass from Matui's criterion to the quotient graph, let $M$ be the
vertex adjacency matrix and let $r_v$ be the number of length-$L$
words forcing $v$. From every initial vertex exactly $r_v$ paths of
length $L$ end at $v$. Hence every row of $M^L$ equals $(r_v)_{v\in G^0}$.
Section~\ref{sec:prefix} uses this identity modulo $d-1$ to relate
equality of the residues $r_v$ to the column sums of $M$.

The same column sums have a dynamical interpretation. For $v\in G^0$, let $V_v\subset X_G$ consist of the paths
starting at $v$, and let $b(v)$ denote the indegree of $v$, counting
edges with multiplicity. Every $x\in V_v$ has precisely $b(v)$
one-step preimages. Such preimage counts also occur in the study of
left special elements for surjective local homeomorphisms
\cite[Section~2.5]{weihe}. Here their residues modulo $d-1$ can be
read in the graph algebra. Indeed, for $s=|G^0|$ we obtain an explicit
isomorphism $\OG\cong M_s(\Od)$, where $\Od$ is the Cuntz
algebra on $d$ generators \cite{cuntz1977}. Under this isomorphism all
vertex projections have the same $K_0$-class $e$ of order $d-1$
\cite{cuntz} (Corollary~\ref{cor:fold-scaled-k}). Write
$c_X(U)=[j_X(1_U)]$ for a clopen set $U\subset X$. Since
$\sigma_G^{-1}(V_v)$ consists of the $b(v)$ edge cylinders ending at
$v$, one has
\[
 c_{X_G}(\sigma_G^{-1}V_v)-c_{X_G}(V_v)=(b(v)-1)e.
\]
On the full shift $c_{\fullshift}(\sigma_d^{-1}U)=c_{\fullshift}(U)$
for every clopen $U$. The homeomorphisms in approximate $K$-conjugacy
transfer this equality to $X_G$. Conversely, when the differences
above vanish, the two factor maps supplied by the folding give
cylinder partitions with the required $K_0$-classes. Shimomura's approximation lemma
\cite[Lemma~3.2]{shimomura} gives homeomorphisms realizing the
partitions. Equivalence of the corresponding nonzero projections in
the purely infinite simple algebras gives the norm approximations on
the diagonals (Corollary~\ref{cor:real-factor}).

\medskip\noindent\textbf{Question 2.}
Which foldings satisfy the finite coefficient $K$-condition with
$\fullshift$ when the homeomorphisms are required only to make the two diagrams
in Definition~\ref{def:finitecoeff} commute on each prescribed pair
of finitely generated subgroups for all sufficiently large indices?

\begin{mainB}
Let $G$ be as in Theorem~A and put $s=|G^0|$. The finite coefficient
$K$-condition for $X_G$ and $\fullshift$ holds if and only if
$\gcd(s,d-1)=1$. This is also equivalent to continuous orbit equivalence
with $\fullshift$ and to a unital isomorphism $\OG\cong\Od$.
\end{mainB}

Indeed, $[1_{\OG}]=se$ generates $K_0(\OG)$ precisely when
$\gcd(s,d-1)=1$. Since $\det(I-M)=1-d$ for every folding on $d$
letters, Matsumoto and Matui's theorem supplies a diagonal-preserving
isomorphism in this case. Its restriction to the diagonal gives a
homeomorphism which, together with its inverse, makes both
diagrams commute on the entire groups of continuous integer-valued functions.

The distinction between the two criteria already appears for $d=4$.
Partition the four letters into classes of sizes $1$ and $3$ and
fold the order-one de Bruijn graph by this partition. The quotient has
two vertices, so $\gcd(2,3)=1$. Its indegrees are $2$ and $6$, and
Theorem~A excludes approximate $K$-conjugacy with $\Sigma_4^+$.
Partitioning the same alphabet into two classes of size $2$ gives
indegrees $4$ and $4$, so both conditions hold
(Example~\ref{ex:six-four}). For a comparison of two foldings on the
same alphabet, Theorem~\ref{thm:six-static} replaces coprimality by
equality of the two greatest common divisors.

The remainder of the article is organized as follows.
Section~\ref{sec:pre} describes foldings and the canonical embeddings
and gives the definitions for a pair of one-sided subshifts.
Section~\ref{sec:prefix} computes the homology classes of the sets
$P_v$ from the adjacency matrix. Section~\ref{sec:algebra} identifies
the graph algebra and computes the maps induced on cylinder classes
by the two factor maps. Section~\ref{sec:real} uses these computations
to prove Theorem~A. Section~\ref{sec:examples} proves Theorem~B and
considers edge shifts outside the class of foldings. In particular,
Corollary~\ref{cor:six-two} gives a two-vertex edge shift approximately
$K$-conjugate to $\Sigma_3^+$ whose topological entropy is $\log 2$.

\section*{Declaration of AI Use}

During the preparation of this manuscript, the author used OpenAI primarily to edit
the language, clarify the exposition, and organize the text. The author
has reviewed and verified all mathematical statements, arguments, and
proofs and takes full responsibility for the content of this paper.

\section{De Bruijn foldings and \texorpdfstring{$K$}{K}-theoretic conditions}\label{sec:pre}

\subsection{Foldings of de Bruijn graphs}\label{sub:foldings}

\smallskip\noindent\textit{Words and the de Bruijn graph.}
Fix an \emph{alphabet} $\mathcal{A}=\{0,\ldots,d-1\}$ of $d\geq2$ letters.
Theorems~A and~B concern $d\geq3$. We include the binary alphabet
in the preliminary constructions for use in Example~\ref{ex:six-binary-trace}.
A word $u=u_1\cdots u_m\in \mathcal{A}^m$ has length $|u|=m$. The set $\mathcal{A}^0=\{\varepsilon\}$
consists of the empty word. Write $uw$ for concatenation and, for
$1\leq L\leq m$, put $\operatorname{suf}_L(u)=u_{m-L+1}\cdots u_m$.

A finite \emph{directed graph} $E$ consists of a finite set $E^0$ of
vertices, a finite set $E^1$ of edges, and two maps $s,r:E^1\to E^0$.
An edge $e$ goes from its \emph{source} $s(e)$ to its \emph{range}
$r(e)$. A labelling by $\mathcal{A}$ assigns to each edge $e$ a letter
$\ell(e)\in \mathcal{A}$. We write $i\xrightarrow{a}j$ for an edge with source
$i$, range $j$, and label $a$. An edge can start and end at the same
vertex (a \emph{loop}), and distinct edges can have the same source and
range.

\begin{df}[De Bruijn graph]\label{def:debruijn}
For $L\geq1$, the \emph{de Bruijn graph of order $L$} on $\mathcal{A}$, denoted
by $B(\mathcal{A},L)$ (written $G(d,L)$ in \cite[Section~3.4]{bleak}), has the words in $\mathcal{A}^L$ as its vertices. Its edges are the
pairs $(u,a)\in \mathcal{A}^L\times \mathcal{A}$. The edge $(u,a)$ is given by
\[
 u_1\cdots u_L \xrightarrow{\ a\ } u_2\cdots u_La
      =\operatorname{suf}_L(ua).
\]
Thus a step deletes the first letter and appends the edge label on the
right. For $L=1$ the new vertex is the one-letter word $a$.
\end{df}

For $\mathcal{A}=\{0,1,2\}$ and $L=2$, there are nine vertices and
$12\xrightarrow{a}2a$ for each $a\in \mathcal{A}$. From any vertex $xy$ we have
$xy\xrightarrow{1}y1\xrightarrow{0}10$. More generally, after reading
$u_1\cdots u_k$ from $x_1\cdots x_L$ the vertex is
$x_{k+1}\cdots x_Lu_1\cdots u_k$ ($0\leq k\leq L$), with empty
strings omitted. At $k=L$ the endpoint is $u_1\cdots u_L$.

\smallskip\noindent\textit{Reading a word in a labelled graph.}
For $m\geq1$, a \emph{path of length $m$} in $E$ is a sequence of edges
$e_1,\ldots,e_m$ with $r(e_k)=s(e_{k+1})$ for $1\leq k<m$.
It starts at $s(e_1)$, ends at $r(e_m)$, and has label word
$\ell(e_1)\cdots\ell(e_m)$. At each vertex we also allow a path of
length zero which stays at that vertex and has label $\varepsilon$.
Suppose each vertex of $E$ has exactly one outgoing edge of each
label in $\mathcal{A}$, as happens in $B(\mathcal{A},L)$. To \emph{read} a word $u$ from
a vertex $i$ means to follow its letters in order, taking at each step
the unique edge carrying the next letter. Write $\delta(i,u)$ for the
vertex reached. In particular, $\delta(i,\varepsilon)=i$, and splitting
a path after the letters of $u$ gives
\begin{equation}\label{eq:read-concatenation}
 \delta(i,uw)=\delta(\delta(i,u),w)
 \qquad\text{for all words }u,w .
\end{equation}

\smallskip\noindent\textit{Identifying vertices.}
A partition of $\mathcal{A}^L$ divides it into nonempty, disjoint subsets whose
union is $\mathcal{A}^L$. Write $u\sim v$ if $u$ and $v$ lie in the same subset.
This is the equivalence relation associated to the partition. Denote the
subset containing $u$ by $\langle u\rangle_\sim$, its \emph{class}.

\begin{df}[Folding]\label{def:debruijn-fold}
Suppose a partition of $\mathcal{A}^L$ satisfies
\begin{equation}\label{eq:folding-relation}
 u\sim v \ \Longrightarrow\
 \operatorname{suf}_L(ua)\sim\operatorname{suf}_L(va)
 \qquad (u,v\in \mathcal{A}^L,\ a\in \mathcal{A}).
\end{equation}
Such an equivalence relation is a \emph{folding} of $B(\mathcal{A},L)$
\cite[Section~3.4]{bleak}. Its \emph{folded automaton} is the labelled quotient
graph whose vertices are the classes and which has the edges
\[
 \langle u\rangle_\sim \xrightarrow{\ a\ }
     \langle\operatorname{suf}_L(ua)\rangle_\sim
 \qquad (a\in \mathcal{A}).
\]
There is one edge for each pair consisting of a source class and a
letter. Condition~\eqref{eq:folding-relation} says precisely that its
range is independent of the representative $u$ of the source class.
We write $G$ for this quotient. In the statements of Theorems~A and~B,
the graph associated with a folding is understood in this sense.
\end{df}

For the nine-vertex graph above, consider the three classes
\[
 A=\{01,02,11,12,21,22\},\qquad B=\{10,20\},\qquad C=\{00\}.
\]
The successor class depends only on the source class and the new letter,
as required by \eqref{eq:folding-relation}.

\[
\begin{array}{c|ccc}
 \text{source}\,\backslash\,\text{label}&0&1&2\\\hline
 A&B&A&A\\
 B&C&A&A\\
 C&C&A&A
\end{array}
\]
In particular, the edges labelled $1$ and $2$ leaving $B$ are distinct
edges from $B$ to $A$. Reading $10$ from any of $A,B,C$ ends at $B$,
whereas reading $00$ from any of them ends at $C$. Reading the single
letter $0$ still depends on the starting vertex.

\smallskip\noindent\textit{Recognizing a folding from its transitions.}
Let $G$ be a finite labelled graph with a nonempty vertex set and
exactly one outgoing edge of each label in $\mathcal{A}$ at each vertex.
In the terminology of \cite[Section~3.3]{bleak}, $G$ is
\emph{synchronizing at level $L$} if there is a map $t:\mathcal{A}^L\to G^0$
such that
\[
 \delta(i,u)=t(u)\qquad (i\in G^0,\ u\in \mathcal{A}^L).
\]
Following \cite[Section~3.3]{bleak}, $t$ is the \emph{synchronizing map}
at level $L$, and $t(u)$ is the state forced by $u$. The automaton is \emph{core} if
$t(\mathcal{A}^L)=G^0$.

\begin{lem}[{\cite[Proposition~3.4]{bleak}}]\label{lem:folding-form}
Let $G$ be as above and fix $L\geq1$. The following are equivalent.
\begin{enumerate}[\rm(i)]
\item\label{fold:quot} $G$ is isomorphic, as a labelled graph, to the quotient by a folding of $B(\mathcal{A},L)$.
\item\label{fold:sync} $G$ is synchronizing at level $L$ and its
synchronizing map is onto.
\end{enumerate}
Here an isomorphism of labelled graphs is a bijection on vertices
together with a bijection on edges preserving sources, ranges, and labels.
\end{lem}

We record the transition formula in our notation. For $u=u_1\cdots u_L$,
$a\in\mathcal A$ and $i\in G^0$, \eqref{eq:read-concatenation} gives
\begin{equation}\label{eq:suffix-transition}
 \delta(t(u),a)=\delta(i,ua)
 =\delta(\delta(i,u_1),\operatorname{suf}_L(ua))
 =t(\operatorname{suf}_L(ua)).
\end{equation}
The last equality uses synchronization of the length-$L$ suffix.

For each folding we fix one such level $L$. All later uses of $t$
refer to the synchronizing map at this chosen level.

\subsection{Cylinder sets in the edge shift and the full shift}\label{sub:blocks}

For the folded automaton $G$ of Subsection~\ref{sub:foldings}, write
$(X,\sigma_X)=(X_G,\sigma_G)$ for its one-sided edge shift and $(\fullshift,\sigma_d)$ for the
one-sided full shift on $\mathcal{A}$.

\smallskip\noindent\textit{The full shift and its cylinders.}
Write $\mathbb N_0=\{0,1,2,\ldots\}$. The \emph{one-sided full shift}
on $\mathcal{A}$ is the space $\fullshift=\mathcal{A}^{\mathbb N_0}$ of all sequences
$z=z_0z_1z_2\cdots$ with $z_j\in \mathcal{A}$, equipped with the map
$\sigma_d(z_0z_1z_2\cdots)=z_1z_2\cdots$. For a finite word $u$ and
$z\in \fullshift$, write $uz$ for the infinite sequence obtained by
placing $u$ before $z$. The \emph{word cylinder} determined by $u$ is
\[
 [u]=\{uz:z\in \fullshift\}.
\]
In particular $[\varepsilon]=\fullshift$. Give $\mathcal{A}$ the discrete
topology and $\fullshift$ the product topology. The word cylinders form
a basis for this topology. For each $m\geq1$ the length-$m$ cylinders
partition $\fullshift$, so each cylinder is both open and closed
(\emph{clopen}). In particular,
\begin{equation}\label{eq:cylinder-refinement}
 [u]=\bigsqcup_{a\in \mathcal{A}}[ua],
\end{equation}
where $\bigsqcup$ denotes disjoint union.

\smallskip\noindent\textit{The edge shift and its coordinates.}
The \emph{one-sided edge shift} of $G$ is the space
\[
 X=\{\,e_0e_1e_2\cdots\in(G^1)^{\mathbb N_0}:
             r(e_j)=s(e_{j+1})\text{ for every }j\geq0\,\}
\]
with the product subspace topology, taking $G^1$ discrete, as in \cite[Section~2.2]{bcw}. The shift deletes the first edge, so $\sigma_X(e_0e_1e_2\cdots)=e_1e_2\cdots$. For a finite path
$\mu$ of positive length, $[\mu]\subset X$ denotes the
\emph{path cylinder} of infinite paths beginning with those particular
edges. A path of length zero at $v$ determines the cylinder of all
paths starting at $v$. Path cylinders are clopen and form a basis,
since a finite initial edge segment specifies a basic neighbourhood.

\begin{lem}\label{lem:path-coordinates}
Giving $G^0$ the discrete topology, the map
\[
 \Theta:X\longrightarrow G^0\times \fullshift,\qquad
 e_0e_1e_2\cdots\longmapsto
  \bigl(s(e_0),\,\ell(e_0)\ell(e_1)\ell(e_2)\cdots\bigr)
\]
is a homeomorphism. In these coordinates the shift is given by
\begin{equation}\label{eq:shift-coordinates}
 (\Theta\sigma_X\Theta^{-1})(i,az)=(\delta(i,a),z)
 \qquad(i\in G^0,\ a\in \mathcal{A},\ z\in \fullshift).
\end{equation}
\end{lem}
\begin{proof}
Each vertex emits exactly one edge of each label, so reading $y$ from
$i$ gives the unique path $\Theta^{-1}(i,y)$. For the path $\mu$ read
from $i$ with word $u$,
\[
 \Theta([\mu])=\{i\}\times[u].
\]
These cylinders form bases on the two sides. Deleting the first edge
of $\Theta^{-1}(i,az)$ gives \eqref{eq:shift-coordinates}.

\end{proof}

The product $\fullshift$ is compact metrizable with a clopen basis.
By \eqref{eq:cylinder-refinement} every cylinder splits into $d\geq2$
nonempty cylinders, so it has no isolated points. Hence $\fullshift$ is
a Cantor space, as is $X\cong G^0\times \fullshift$ by
Lemma~\ref{lem:path-coordinates}.

\smallskip\noindent\textit{Clopen sets determined by vertices.}
For $v\in G^0$, define
\begin{equation}\label{eq:two-blocks}
 \begin{aligned}
 V_v&=\{x=e_0e_1\cdots\in X:s(e_0)=v\},\\
 P_v&=\{z\in \fullshift:t(z_0\cdots z_{L-1})=v\}
      =\bigsqcup_{u\in t^{-1}(\{v\})}[u],\qquad
 r_v=|t^{-1}(\{v\})|.
 \end{aligned}
\end{equation}
Since $t$ is onto, \eqref{eq:two-blocks} gives the clopen partitions
\[
 X=\bigsqcup_{v\in G^0}V_v,\qquad \fullshift=\bigsqcup_{v\in G^0}P_v,
 \qquad r_v\geq1,\quad\sum_v r_v=d^L.
\]
Thus $P_v$ is the union of the word cylinders whose initial length-$L$
word forces the state $v$. Here $|\cdot|$ denotes cardinality and $\Theta(V_v)=\{v\}\times \fullshift$.
For $u\in \mathcal{A}^L$, iterating \eqref{eq:shift-coordinates} gives
$(\Theta\sigma_X^L\Theta^{-1})(i,uz)=(t(u),z)$. Thus
\begin{equation}\label{eq:terminal-vertex-preimage}
 \Theta\bigl((\sigma_X^L)^{-1}(V_v)\bigr)=G^0\times P_v .
\end{equation}

For the folding with vertex set $\{A,B,C\}$ in
Subsection~\ref{sub:foldings}, take $L=2$. The three sets are
\[
 \begin{aligned}
 P_A&=[01]\sqcup[02]\sqcup[11]\sqcup[12]\sqcup[21]\sqcup[22],\\
 P_B&=[10]\sqcup[20],\qquad P_C=[00].
 \end{aligned}
\]
For example, a label sequence $10z$ lies in $P_B$. There are three
paths with that label sequence, one from each initial vertex. Only
the path starting at $B$ belongs to $V_B$, and all three are at $B$
after the first two edges.

\begin{rem}[Choice of synchronizing level]\label{rem:terminal-level}
Write $t_L$, $P_v^{(L)}$ and $r_v^{(L)}$ for the objects in
\eqref{eq:two-blocks} when the level needs to be indicated. If $k\geq0$,
$w\in \mathcal{A}^k$ and $u\in \mathcal{A}^L$, reading $wu$ from any vertex ends at
$t_L(u)$, since the final $L$ letters force that endpoint. Thus $L+k$
is also a synchronizing level, its synchronizing map is onto, and
\begin{equation}\label{eq:terminal-level-change}
 t_{L+k}(wu)=t_L(u),\qquad
 P_v^{(L+k)}=(\sigma_d^k)^{-1}(P_v^{(L)}),\qquad
 r_v^{(L+k)}=d^k r_v^{(L)} .
\end{equation}
Indeed, the words ending at $v$ at level $L+k$ are precisely
$wu$ with $w\in \mathcal{A}^k$ and $u\in t_L^{-1}(\{v\})$, by
\eqref{eq:read-concatenation}. For example,
$P_C^{(2)}=[00]$ and $P_C^{(3)}=[000]\sqcup[100]\sqcup[200]$.
Thus the block can change, while
$r_v^{(L+k)}\equiv r_v^{(L)}\pmod{d-1}$ by
\eqref{eq:terminal-level-change}.

\end{rem}

\smallskip\noindent\textit{Counting paths.}
Following \cite[Definition~2.2.3]{lindmarcus}, the adjacency matrix $M=(M_{ij})_{i,j\in G^0}$ has entry $M_{ij}$
equal to the number of edges from $i$ to $j$. Thus
$b(v)=\sum_{i\in G^0}M_{iv}$, the \emph{indegree} of $v$, counts every
edge ending at $v$, including each of several edges with the same
endpoints. Write $\mathbf1$ for the column vector with all entries $1$
indexed by $G^0$, and $\mathbf r=(r_v)_{v\in G^0}$ for the row of
cardinalities $|t^{-1}(\{v\})|$.

\begin{lem}\label{lem:fold-count}
The adjacency matrix and the one-step preimages satisfy
\begin{equation}\label{eq:folding-counts}
 M\mathbf1=d\mathbf1,\qquad M^L=\mathbf1\mathbf r,\qquad
 |\sigma_X^{-1}(\{x\})|=b(v)\quad(x\in V_v).
\end{equation}
\end{lem}
\begin{proof}
At $i$ there is one outgoing edge per letter, so $\sum_j M_{ij}=d$.
Matrix multiplication counts successive edges. At length $L$ the
label map on paths starting at $i$ therefore gives
\[
 (M^L)_{iv}
 =\#\{\mu:|\mu|=L,\ s(\mu)=i,\ r(\mu)=v\}
 =\#\{u\in \mathcal{A}^L:t(u)=v\}=r_v .
\]
For $x\in V_v$, $e\mapsto ex$ is a bijection from $r^{-1}(\{v\})$
onto $\sigma_X^{-1}(\{x\})$. This proves \eqref{eq:folding-counts}.
\end{proof}

More generally, $\sigma_X|_{[e]}:[e]\to V_{r(e)}$ has continuous
inverse $x\mapsto ex$ \cite[Section~2.2]{bcw}. Every $v$ is the endpoint
of a length-$L$ path by \eqref{eq:two-blocks}, so $b(v)>0$ and
$\sigma_X$ is a surjective local homeomorphism.

In the running example, ordering the vertices as $A,B,C$ gives
\[
 M=\begin{pmatrix}2&1&0\\2&0&1\\2&0&1\end{pmatrix},
 \qquad M^2=\begin{pmatrix}6&2&1\\6&2&1\\6&2&1\end{pmatrix}.
\]
The column sums of $M$ give the one-step preimage counts $(6,1,2)$.
Each row of $M^2$ gives the numbers of length-two words forcing each vertex $(6,2,1)$.

\smallskip\noindent\textit{Finite prefix replacements.}
A finite word $u$ is a \emph{prefix} of a word or sequence if it
occupies its first $|u|$ positions. A \emph{finite complete prefix
code} is a finite set $\mathcal U$ of words such that every infinite
sequence begins with exactly one member of $\mathcal U$. Its cylinders
therefore partition $\fullshift$. Given two such codes $\mathcal U$
and $\mathcal W$ and a bijection $\theta:\mathcal U\to\mathcal W$, set
\[
 h(uz)=\theta(u)z\qquad(u\in\mathcal U,\ z\in \fullshift).
\]
On $[u]$ this rule replaces the displayed prefix and keeps the tail
$z$. Its inverse uses $\theta^{-1}$. Since the cylinders on each
side form a finite clopen partition, $h$ is a homeomorphism. The
group of all such homeomorphisms is the Higman--Thompson group
$V_{d,1}$ in its full-shift realization \cite[Remark~6.3 and Section~6.7.1]{matuifull}.
Subdividing a matched pair $u,\theta(u)$ into the pairs
$ua,\theta(u)a$ for $a\in \mathcal{A}$ gives the same homeomorphism.

We use Matui's term \emph{transposition} \cite[Section~6.4]{matuifull}.
For $v\ne w$, a transposition interchanging $P_v$ and $P_w$ is an
$h\in V_{d,1}$ satisfying
\[
 h^2=\id_{\fullshift},\qquad h(P_v)=P_w,\qquad
 h(y)=y\ \ (y\in \fullshift\setminus(P_v\cup P_w)).
\]
In particular $h(P_w)=P_v$.

For example, refine the two cylinders of $P_B$ in the running example
by \eqref{eq:cylinder-refinement}. Their six new prefixes are $i0a$,
where $i\in\{1,2\}$ and $a\in \mathcal{A}$. The six prefixes of $P_A$ are
precisely the words $ai$ with the same indices. Exchange each pair by
\[
 i0az\longleftrightarrow aiz
 \qquad (i\in\{1,2\},\ a\in \mathcal{A},\ z\in \fullshift)
\]
and fix $[00]=P_C$ pointwise. The twelve indicated cylinders together
with $[00]$ form a partition, so the rule defines a transposition interchanging $P_A$ and $P_B$. Section~\ref{sec:prefix} gives the criterion for
such an exchange for any two distinct sets $P_v,P_w$.

\subsection{Diagonal embeddings and cylinder classes}\label{sub:diagonal}

Let $E_d$ be the graph with one vertex and one loop for each
$a\in\mathcal{A}$. Its edge paths are identified with $\fullshift$ by their
labels. The following construction applies to the two pairs
$(G,X)$ and $(E_d,\fullshift)$. We state the formulas for $(G,X)$. Replacing
$(G,X)$ by $(E_d,\fullshift)$ gives the corresponding formulas for the full shift.

\smallskip\noindent\textit{The graph algebra.}
Write $\mathcal O_X=C^*(G)$ for the graph $C^*$-algebra, equivalently
the Cuntz--Pimsner algebra of the graph correspondence \cite[Examples~8.3, 8.7 and 8.13]{raeburn}. Here $s(e)$ is
the initial vertex and $r(e)$ the terminal vertex. In translating Raeburn's graph correspondence to our convention one interchanges $r$ and $s$. The correspondence is $C(G^1)$ over $C(G^0)$, with
actions and inner product given, for $f,g\in C(G^0)$ and
$\xi,\zeta\in C(G^1)$, by
\[
 (f\cdot\xi\cdot g)(e)=f(s(e))\xi(e)g(r(e)),\qquad
 \langle\xi,\zeta\rangle(v)=\sum_{r(e)=v}\overline{\xi(e)}\zeta(e).
\]
We use the universal presentation of its algebra with mutually orthogonal
vertex projections $p_v$ ($v\in G^0$) and edge partial isometries
$S_e$ ($e\in G^1$) whose range projections are mutually orthogonal,
subject to
\begin{equation}\label{eq:graph-relations}
 S_e^*S_e=p_{r(e)},\qquad
 p_v=\sum_{s(e)=v}S_eS_e^*.
\end{equation}
A projection satisfies $p=p^*=p^2$. For a partial isometry $S$, the
projections $S^*S$ and $SS^*$ are its initial and range projections.
Every vertex here emits $d$ edges, so the second relation applies at
every vertex, and $1_{\mathcal O_X}=\sum_{v\in G^0}p_v$.
In particular $\mathcal O_X=\OG$ and $C^*(E_d)=\Od$. In the one-vertex case the vertex projection
is the unit and the loop generators are $d$ isometries with orthogonal
ranges summing to the unit. These are the defining relations of $\Od$ \cite{cuntz1977}.

For a path $\mu=e_1\cdots e_m$ of positive length, set
$s(\mu)=s(e_1)$, $r(\mu)=r(e_m)$ and $S_\mu=S_{e_1}\cdots S_{e_m}$.
For a path of length zero at $v$, put $s(v)=r(v)=v$ and $S_v=p_v$.
Let $q_\mu=S_\mu S_\mu^*$. We write $C(X)$ for the algebra of
continuous complex-valued functions on $X$ with the supremum norm
and $1_U$ for the indicator function of a clopen set $U\subset X$.

\begin{lem}[{\cite[Theorem~3.7]{webster}}]\label{lem:diagonal-embedding}
There is a unique unital injective $*$-homomorphism
$j_X:C(X)\to\mathcal O_X$ satisfying $j_X(1_{[\mu]})=q_\mu$ for every
finite path $\mu$. Its range is the closed linear span of the $q_\mu$,
called the canonical diagonal of $\mathcal O_X$.
\end{lem}
Each vertex emits a finite nonzero number of edges. Thus Webster's
boundary path space is $X$ for $G$ and $\fullshift$ for $E_d$. In our path convention the
Gelfand transform on the diagonal sends $q_\mu$ to $1_{[\mu]}$.
Its inverse followed by inclusion gives $j_X$, and the vertex
projections sum to the unit. Applying $j_X$ to the partition of $[\mu]$ by its next edge gives
\begin{equation}\label{eq:path-refinement}
 q_\mu=\sum_{s(e)=r(\mu)}q_{\mu e}.
\end{equation}

\smallskip\noindent\textit{Classes of clopen sets.}
For a unital $C^*$-algebra $A$, write $K_0(A)$ for the Grothendieck
group of stable Murray--von Neumann classes of projections in matrix
algebras over $A$. Here projections are equivalent if they are the
initial and range projections of a partial isometry, after adjoining
zero blocks as needed. Write $[p]$ for the class of $p$. Addition is
induced by block sum, so $[p+q]=[p]+[q]$ for orthogonal projections.
We use the ordinary projection and unitary descriptions of $K$-theory \cite{blackadar}. Write $K_1(A)$ for the group of stable homotopy classes of unitaries,
and $\varphi_{*i}$ for the map on $K_i$ induced by a unital
$*$-homomorphism $\varphi$, $i=0,1$. For a clopen $U\subset X$ set
\[
 c_X(U)=[j_X(1_U)]\in K_0(\mathcal O_X).
\]

\begin{lem}\label{lem:clopen-classes}
Pointwise rank identifies $K_0(C(X))$ with $C(X,\Z)$ and sends
$[1_U]$ to $1_U$ for every clopen $U\subset X$. Also $K_1(C(X))=0$.
For every finite path $\mu$ and every finite disjoint cylinder
decomposition $U=\bigsqcup_{j=1}^{\ell}[\mu_j]$, one has
\begin{equation}\label{eq:path-class}
 c_X([\mu])=[p_{r(\mu)}],\qquad
 c_X(U)=\sum_{j=1}^{\ell}[p_{r(\mu_j)}].
\end{equation}
\end{lem}
\begin{proof}
Since $X$ is a Cantor space, the assertions about $K_*(C(X))$ follow
from the calculation recalled in \cite[Section~1, preceding Theorem~1.1]{putnam}.
The pointwise rank of $1_U$ is one on $U$ and zero elsewhere.

The path relations \cite[Corollary~1.14 and the preceding path formula]{raeburn},
in our source and range convention, and Lemma~\ref{lem:diagonal-embedding} give
\[
  S_\mu^*S_\mu=p_{r(\mu)},\qquad
  S_\mu S_\mu^*=j_X(1_{[\mu]}).
\]
The two projections are Murray--von Neumann equivalent via $S_\mu$.
Every nonempty clopen $U$ has a finite cover by cylinders contained in $U$
by compactness. Refining to their maximum length gives a disjoint
cylinder decomposition. For $U=\varnothing$ take the empty union. For disjoint
cylinders, $1_U=\sum_j1_{[\mu_j]}$ is an orthogonal sum, proving
\eqref{eq:path-class}.
\end{proof}

For $\fullshift$ put $\bar1=[1_{\Od}]$. Applying \eqref{eq:path-class} to the
word cylinders in \eqref{eq:two-blocks} and to the length-zero path at
$v$ gives, respectively,
\[
 c_{\fullshift}(P_v)=r_v\bar1,\qquad c_X(V_v)=[p_v].
\]

\subsection{Approximate \texorpdfstring{$K$}{K}-conjugacy and the finite coefficient
\texorpdfstring{$K$}{K}-condition}\label{sub:compatibility}

In this subsection we consider two nonempty one-sided subshifts $(X,\sigma_X)$ and $(Y,\sigma_Y)$
with compatible metrics $d_X,d_Y$. Here a subshift over a finite alphabet
$\mathcal A$ is a closed subset $Z\subseteq\mathcal A^{\mathbb N_0}$
satisfying $\sigma(Z)\subseteq Z$, with $\sigma_Z=\sigma|_Z$
\cite[Section~1.1]{brixcovers}. The two alphabets may differ.

For $Z=X,Y$, let $\widetilde Z$ be its cover and
$\pi_Z:\widetilde Z\to Z$ its canonical factor map in the sense of
\cite[Definition~2.1 and Section~2.1]{brixcovers}.
In the associated shift algebra $\mathcal O_Z$ the diagonal is
$\mathcal D_Z=C(\widetilde Z)$, and the coefficient embedding is
\[
 j_Z:C(Z)\longrightarrow\mathcal D_Z\subseteq\mathcal O_Z,\qquad
 j_Z(f)=f\circ\pi_Z
\]
\cite[Sections~2.1 and~2.3]{brixcovers}. For the edge shifts
of Subsection~\ref{sub:diagonal}, $\pi_Z$ is a homeomorphism
\cite[Section~2.1]{brixcovers}. The shift groupoid is then the graph
groupoid, and the canonical isomorphism from its $C^*$-algebra to the graph algebra sends
the cylinder projections to the projections of
Lemma~\ref{lem:diagonal-embedding} \cite[Proposition~2.2]{bcw}.
For a clopen $U\subset X$ put $c_X(U)=[j_X(1_U)]$, and define $c_Y$
similarly. All coefficient conditions below concern these embeddings.

For $Z=X,Y$, the nonempty length-$m$ cylinders form a finite clopen
partition $\mathcal P_m(Z)$. These partitions refine and generate the
topology. The spans of their indicators therefore have dense union in
$C(Z)$. By the finite-dimensional calculation and continuity of
$K$-theory \cite[5.2.4, 5.5.1 and 8.1.5]{blackadar}, pointwise rank gives
\[
 K_0(C(Z))=C(Z,\Z),\qquad [1_U]=1_U,\qquad K_1(C(Z))=0
 \quad(U\subset Z\text{ clopen}).
\]
Indeed, the map on integer-valued functions under refinement repeats
each value on the subcylinders, and every continuous integer-valued
function is constant on cylinders of some common length.
For a homeomorphism $h:X\to Y$, write $h^*(f)=f\circ h$.
In degree zero $(h^*)_{*0}(f)=f\circ h$ for $f\in C(Y,\Z)$.
In particular, $1_U\circ h=1_{h^{-1}(U)}$ for clopen $U\subset Y$.

\smallskip\noindent\textit{Approximate $K$-conjugacy.}
The shift algebras $\mathcal O_X,\mathcal O_Y$ are separable, nuclear
and satisfy the universal coefficient theorem (UCT)
\cite[Theorem~17 in the preprint version]{carlsensilvestrov}.
For separable $C^*$-algebras $A,B$ with $A$ in the bootstrap class,
let $\iota_{\mathrm{UCT}}$ denote the injective homomorphism from the degree-zero
$\operatorname{Ext}^1$ term into $KK(A,B)$ in the Rosenberg--Schochet
sequence \cite[Theorem~1.17]{rs}. We use the quotient
\cite[Section~1.6]{dadarlat}
\[
 KL(A,B)=KK(A,B)\big/
 \iota_{\mathrm{UCT}}\!\left(
 \bigoplus_{i=0,1}\operatorname{Pext}^1_{\Z}(K_i(A),K_{1-i}(B))\right).
\]
A pure extension splits over every finitely generated
subgroup of its quotient. Write $[\Phi]_{KL}$ for the class of a
$*$-homomorphism $\Phi:A\to B$. If $[\Phi]_{KL}=\kappa$, then
$\Phi_{*i}=\kappa_{*i}$ for $i=0,1$. For invertible $\kappa$,
$\kappa^{-1}\in KL(B,A)$ denotes its inverse under the Kasparov product.

\medskip

Lin and Matui introduced approximate $K$-conjugacy for minimal Cantor
homeomorphisms \cite[Definition~5.3]{linmatui}.
\begin{df}\label{def:ak}
Following the formulation for Cantor--circle extensions in
\cite[Section~2]{llw}, we say that $(X,\sigma_X)$ and $(Y,\sigma_Y)$ are
\emph{approximately $K$-conjugate} if there exist homeomorphisms
$\phi_n:X\to Y$ and $\psi_n:Y\to X$, unital $*$-isomorphisms
$\Phi_n:\mathcal O_X\to\mathcal O_Y$ and $\Psi_n:\mathcal O_Y\to\mathcal O_X$
for $n\geq1$, and a single
invertible class $\kappa\in KL(\mathcal O_Y,\mathcal O_X)$ with
$[\Psi_n]_{KL}=\kappa$ and $[\Phi_n]_{KL}=\kappa^{-1}$ for all $n$,
having the following properties.
\begin{enumerate}[\rm(i)]
\item The homeomorphisms approximately intertwine the shifts in both
directions:
\begin{equation}\label{eq:ak-uniform}
 \begin{aligned}
 \sup_{x\in X}d_{Y}(\phi_n(\sigma_X x),\sigma_Y(\phi_n x))&\longrightarrow0,\\
 \sup_{y\in Y}d_X(\psi_n(\sigma_Y y),\sigma_X(\psi_n y))&\longrightarrow0 .
 \end{aligned}
\end{equation}
\item On the canonical embeddings $j_Y$ and $j_X$ one has pointwise norm convergence:
\begin{equation}\label{eq:ak-norm}
 \begin{aligned}
 \|\Psi_n(j_{Y}(f))-j_X(f\circ \phi_n)\|_{\mathcal O_X}
    &\longrightarrow0 &&(f\in C(Y)),\\
 \|\Phi_n(j_X(g))-j_{Y}(g\circ \psi_n)\|_{\mathcal O_Y}
    &\longrightarrow0 &&(g\in C(X)).
 \end{aligned}
\end{equation}
\end{enumerate}
\end{df}

By the substitutions $y=\phi_n(x)$ and $x=\psi_n(y)$,
\eqref{eq:ak-uniform} is uniform convergence of
$\phi_n\sigma_X\phi_n^{-1}$ to $\sigma_Y$ and of $\psi_n\sigma_Y\psi_n^{-1}$ to
$\sigma_X$. Compactness makes these uniform limits independent of the
compatible metrics. Definition~\ref{def:ak} uses the same homeomorphisms
in \eqref{eq:ak-uniform} and \eqref{eq:ak-norm}. The algebra isomorphisms have mutually inverse $KL$-classes.
The definition allows $\psi_n$ and $\phi_n^{-1}$ to differ.

\smallskip\noindent\textit{Finite coefficient $K$-condition.}
Put $K_*(A)=K_0(A)\oplus K_1(A)$, graded by the two summands.
A unit-preserving graded isomorphism
$\eta:K_*(\mathcal O_Y)\to K_*(\mathcal O_X)$ consists of group
isomorphisms in each degree with
$\eta_0([1_{\mathcal O_Y}])=[1_{\mathcal O_X}]$. Write $k=\eta_0$.
For finitely generated subgroups $H\subset C(Y,\Z)$ and
$J\subset C(X,\Z)$, the two diagrams are
\[
\begin{CD}
 H @>{\phi_n^*}>> C(X,\Z)\\
 @V{(j_Y)_{*0}}VV @VV{(j_X)_{*0}}V\\
 K_0(\mathcal O_Y) @>{k}>> K_0(\mathcal O_X)
\end{CD}
\qquad\text{and}\qquad 
\begin{CD}
 J @>{\psi_n^*}>> C(Y,\Z)\\
 @V{(j_X)_{*0}}VV @VV{(j_Y)_{*0}}V\\
 K_0(\mathcal O_X) @>{k^{-1}}>> K_0(\mathcal O_Y).
\end{CD}
\]
Each arrow out of $H$ or $J$ denotes the restriction to that subgroup. The upper arrows act by composition.
For Cantor--circle extensions, the finite coefficient $K$-condition
in \cite[Section~2, main theorem, condition~\textup{(2)}]{llw}
uses an ordered, unit-preserving graded isomorphism and involves
both coefficient degrees. In the definition below $\eta$ is required
to preserve the grading and the unit class. The coefficient diagrams
occur only in degree zero because $K_1(C(X))=K_1(C(Y))=0$.
The definition also requires $\eta_1:K_1(\mathcal O_Y)\to K_1(\mathcal O_X)$
to be an isomorphism.

\begin{df}\label{def:finitecoeff}
The shifts $(X,\sigma_X)$ and $(Y,\sigma_Y)$ satisfy the \emph{finite coefficient $K$-condition} if
there exist a unit-preserving graded isomorphism $\eta$ as above and
sequences of homeomorphisms $\phi_n:X\to Y$, $\psi_n:Y\to X$ such that for
every pair of finitely generated subgroups $H\subset C(Y,\Z)$,
$J\subset C(X,\Z)$ there is $N=N(H,J)$ with the following equalities
for every $n\geq N$:
\begin{align}
 (j_X)_{*0}(f\circ \phi_n)&=k((j_Y)_{*0}(f)) &&(f\in H),\label{eq:finite-forward}\\
 (j_Y)_{*0}(a\circ \psi_n)&=k^{-1}((j_X)_{*0}(a)) &&(a\in J).\label{eq:finite-reverse}
\end{align}
\end{df}

The subgroups in this definition are taken in the integral $K_0$ of
the coefficient algebras $C(Y)$ and $C(X)$. The index $N$ works for both diagrams and for
every element of the two prescribed subgroups.

\begin{lem}\label{lem:atoms}
Fix $\eta,\phi_n,\psi_n$ with the domains and codomains of
Definition~\ref{def:finitecoeff}. The following are equivalent.
\begin{enumerate}[\rm(i)]
\item\label{fc:groups} Equations~\eqref{eq:finite-forward}
and~\eqref{eq:finite-reverse} hold eventually on every prescribed
pair $H,J$ of finitely generated subgroups.
\item\label{fc:atoms} For every finite clopen partition $\mathcal U$
of $Y$ and $\mathcal V$ of $X$, there is an index $N$ such that for
all $n\geq N$,
\begin{equation}\label{eq:atom-diagrams}
 \begin{aligned}
 c_X(\phi_n^{-1}(U))&=k(c_Y(U)) &&(U\in\mathcal U),\\
 c_Y(\psi_n^{-1}(V))&=k^{-1}(c_X(V)) &&(V\in\mathcal V).
 \end{aligned}
\end{equation}
\end{enumerate}
\end{lem}
\begin{proof}
\proofstep{fc:groups}{fc:atoms}
Take $H=\langle1_U:U\in\mathcal U\rangle$ and
$J=\langle1_V:V\in\mathcal V\rangle$. Applying
\eqref{eq:finite-forward} and \eqref{eq:finite-reverse} to these
generators gives \eqref{eq:atom-diagrams}, since
$1_U\circ \phi_n=1_{\phi_n^{-1}(U)}$ and $1_V\circ \psi_n=1_{\psi_n^{-1}(V)}$.

\proofstep{fc:atoms}{fc:groups}
On each space take the common refinement of the clopen level sets
of a finite generating set for the respective subgroup. These partitions
are finite because each generator is continuous into the discrete group
$\Z$ and has compact image. For the zero
subgroup take the one-atom partition. Each $f\in H$ and $a\in J$ then
has the form $f=\sum_U m_U1_U$, $a=\sum_V l_V1_V$, with integer
coefficients. For all $n\geq N(\mathcal U,\mathcal V)$,
\eqref{eq:atom-diagrams} and additivity give
\[
 \begin{aligned}
 (j_X)_{*0}(f\circ \phi_n)
 &=\sum_U m_U c_X(\phi_n^{-1}(U))
   =k\!\left(\sum_U m_U c_Y(U)\right)=k((j_Y)_{*0}(f)),\\
 (j_Y)_{*0}(a\circ \psi_n)
 &=\sum_V l_V c_Y(\psi_n^{-1}(V))
   =k^{-1}\!\left(\sum_V l_V c_X(V)\right)=k^{-1}((j_X)_{*0}(a)).
 \end{aligned}
\]
\end{proof}

\begin{cor}\label{cor:ak-finite}
If $(X,\sigma_X)$ and $(Y,\sigma_Y)$ are approximately $K$-conjugate, the same homeomorphisms
realize the finite coefficient $K$-condition for $\eta=\kappa_*$.
\end{cor}
\begin{proof}
Since $\Psi_n$ and $\Phi_n$ are unital isomorphisms with fixed inverse
$KL$-classes, $\eta=\kappa_*=(\Psi_n)_*$ is a unit-preserving graded
isomorphism with inverse $(\Phi_n)_*$. Put $k=\eta_0$. For the finite partitions $\mathcal U,\mathcal V$ of
Lemma~\ref{lem:atoms}, \eqref{eq:ak-norm} gives a common $N$ such that,
for all $n\geq N$,
\[
 \begin{aligned}
 \max_{U\in\mathcal U}\|\Psi_n(j_Y(1_U))-j_X(1_{\phi_n^{-1}(U)})\|&<1,\\
 \max_{V\in\mathcal V}\|\Phi_n(j_X(1_V))-j_Y(1_{\psi_n^{-1}(V)})\|&<1.
 \end{aligned}
\]
By \cite[Proposition~4.6.6]{blackadar} each pair of projections is
unitarily equivalent. Hence
\[
 \begin{aligned}
 k(c_Y(U))&=[\Psi_n(j_Y(1_U))]
       =[j_X(1_{\phi_n^{-1}(U)})]=c_X(\phi_n^{-1}(U)),\\
 k^{-1}(c_X(V))&=[\Phi_n(j_X(1_V))]
       =[j_Y(1_{\psi_n^{-1}(V)})]=c_Y(\psi_n^{-1}(V)).
 \end{aligned}
\]
These are the two identities in \eqref{eq:atom-diagrams}.
Lemma~\ref{lem:atoms} completes the proof.
\end{proof}

The implication in Corollary~\ref{cor:ak-finite} uses the coefficient norm
limits. The additional requirement in \eqref{eq:ak-uniform} concerns
the action of the shifts on these same homeomorphisms. This is the
compatibility used in Section~\ref{sec:real}.

\section{Transpositions and vertex indegrees}\label{sec:prefix}

Let $G$ be the quotient by a folding of $B(\mathcal A,L)$, where $|\mathcal A|=d\geq3$.
Write $s=|G^0|$ and $R=\Z/(d-1)\Z$. A bar denotes reduction modulo
$d-1$, entrywise for vectors and matrices. The clopen partition
$\{P_v:v\in G^0\}$ and the integers $r_v$ are defined in
\eqref{eq:two-blocks}. For a given pair $P_v,P_w$, Matui's criterion reduces the
existence of a transposition to equality of two classes in the homology
of the full shift. We compute these classes from the adjacency matrix
$M$ of $G$. For transpositions of all pairs, only the column sums of $M$ are needed.

\subsection{Transpositions of clopen sets}\label{sub:prefix-clopen}

We use transpositions as in Subsection~\ref{sub:blocks}, now for any
two disjoint nonempty clopen sets $U,W\subseteq\fullshift$. A transposition interchanging them sends $U$ onto $W$ and has support exactly $U\cup W$.

Write $\mathcal G_d$ for the groupoid of the one-sided full $d$-shift.
For a clopen set $U$, write $[1_U]_{\mathcal G_d}$ for the class
of its indicator in the zeroth homology group $H_0(\mathcal G_d)$.
The groupoid here is that of the full shift, where the sets $P_v$
lie. Matui identifies its topological full group with $V_{d,1}$
\cite[Remark~6.3]{matuifull}. Under the identification
\[
 H_0(\mathcal G_d)\cong R,\qquad [1_{\fullshift}]_{\mathcal G_d}
 \longmapsto \overline{1},
\]
each word cylinder $[u]$ has class $\overline{1}$. Indeed, the
compact open bisection inducing $z\mapsto uz$ has source $\fullshift$
and range $[u]$. In degree zero the difference of the indicators of
the range and source of a compact open bisection is a boundary
\cite[Sections~6.1,~6.2 and~6.7.1]{matuifull}. By additivity, any
finite disjoint cylinder decomposition $U=\bigsqcup_{u\in F}[u]$ gives
\begin{equation}\label{eq:prefix-type}
 [1_U]_{\mathcal G_d}=\overline{|F|}\quad\text{in }R .
\end{equation}

\begin{prop}[{\cite[Remark~6.3, Theorem~6.4 and Section~6.4]{matuifull}}]
\label{prop:prefix-clopen}
Let $U,W\subseteq\fullshift$ be disjoint nonempty clopen sets. There
exists a transposition in $V_{d,1}$ interchanging $U$ and $W$ if and
only if $[1_U]_{\mathcal G_d}=[1_W]_{\mathcal G_d}$.
\end{prop}

\subsection{A criterion for \texorpdfstring{$P_v$ and $P_w$}{a pair of vertices}}\label{sub:prefix-pair}

We now apply Proposition~\ref{prop:prefix-clopen} to the partition
defined by the folding. The set $P_v$ is the disjoint union of the
$r_v$ length-$L$ cylinders in \eqref{eq:two-blocks}, so
\eqref{eq:prefix-type} gives
\begin{equation}\label{eq:prefix-type-block}
 [1_{P_v}]_{\mathcal G_d}=\overline{r_v}\qquad(v\in G^0).
\end{equation}
To compare these classes using only $M$, we use the identity
$M^L=\mathbf1\mathbf r$ from \eqref{eq:folding-counts}. Here
$\mathbf r=(r_v)_{v\in G^0}$ is a row and $\mathbf1$ is the column
of $s$ ones, over $\Z$ or $R$ as appropriate.

\begin{lem}\label{lem:prefix-row}
The system of equations over $R$
\begin{equation}\label{eq:prefix-row}
 a\overline M=a,\qquad a\mathbf1=1
\end{equation}
has a unique row-vector solution $a\in R^{1\times s}$, namely
$a=\overline{\mathbf r}$.
\end{lem}
\begin{proof}
Over $\Z$, \eqref{eq:folding-counts} and $M^LM=MM^L$ give
\[
 \mathbf1(\mathbf rM)=M^LM=MM^L=d\mathbf1\mathbf r.
\]
Comparing any row gives $\mathbf rM=d\mathbf r$.
Also $\mathbf r\mathbf1=d^L$ by \eqref{eq:two-blocks}. As
$d\equiv1\pmod{d-1}$, reduction gives \eqref{eq:prefix-row} for
$\overline{\mathbf r}$. Conversely, every solution satisfies
\[
 a=a\overline M^{\,L}
   =a(\mathbf1\overline{\mathbf r})
   =(a\mathbf1)\overline{\mathbf r}=\overline{\mathbf r}.
\]
\end{proof}

\begin{cor}\label{thm:prefix-pair}
Let $v,w\in G^0$ be distinct and let $a$ be the solution of
\eqref{eq:prefix-row}. The following conditions are equivalent.
\begin{enumerate}[\rm(i)]
\item\label{pair:swap} A transposition in $V_{d,1}$ interchanges $P_v$ and $P_w$.
\item\label{pair:count} $r_v\equiv r_w\pmod{d-1}$.
\item\label{pair:matrix} $a_v=a_w$ in $R$.
\end{enumerate}
\end{cor}

\begin{proof}
\proofstep{pair:swap}{pair:count}
By Proposition~\ref{prop:prefix-clopen}, the transposition gives
$[1_{P_v}]_{\mathcal G_d}=[1_{P_w}]_{\mathcal G_d}$. Apply
\eqref{eq:prefix-type-block}.

\proofstep{pair:count}{pair:matrix}
Lemma~\ref{lem:prefix-row} gives $a_v=\overline{r_v}$ and
$a_w=\overline{r_w}$.

\proofstep{pair:matrix}{pair:swap}
By the same lemma and \eqref{eq:prefix-type-block}, the two sets have
equal homology classes. They are nonempty and disjoint by
\eqref{eq:two-blocks}, so Proposition~\ref{prop:prefix-clopen}
gives the required transposition.
\end{proof}

Fix $i\in G^0$, and let $e_i$ be the row with a $1$ at $i$ and
zeros elsewhere. Put $a^{(j)}=e_i\overline M^{\,j}$ for $j\ge0$.
By \eqref{eq:folding-counts}, $a^{(j)}\mathbf1=1$ for every $j$.
If $a^{(j+1)}=a^{(j)}$, then $a^{(j)}\overline M=a^{(j)}$,
so Lemma~\ref{lem:prefix-row} gives $a^{(j)}=a$.
Such an equality occurs by $j=L$, since
$a^{(L)}=\overline{\mathbf r}$ and
$\overline{\mathbf r}\,\overline M=\overline{\mathbf r}$.
Thus the first two equal consecutive rows determine $a$ without prior
knowledge of $L$.

\subsection{The indegree criterion}\label{sub:prefix-all}

By Corollary~\ref{thm:prefix-pair}, every pair of distinct blocks
can be interchanged precisely when all coordinates of $a$ are equal.
A constant row is a multiple of $\mathbf1^{\mathsf T}$, and the
$v$-entry of $\mathbf1^{\mathsf T}M$ is the indegree $b(v)$.
This gives a criterion directly in terms of the edges ending at each
vertex, counted with multiplicity.

\begin{thm}\label{thm:prefix-all}
The following conditions are equivalent.
\begin{enumerate}[\rm(i)]
\item\label{all:swap} For each pair of distinct vertices $v,w$, there is a transposition in $V_{d,1}$ interchanging $P_v$ and $P_w$.
\item\label{all:uniform} All integers $r_v$, $v\in G^0$, are
congruent modulo $d-1$.
\item\label{all:degree} $b(v)\equiv1\pmod{d-1}$ for every $v\in G^0$.
\end{enumerate}
Whenever these conditions hold, $\gcd(s,d-1)=1$ and
$\overline{r_v}=\overline{s}^{\,-1}$ for every $v$.
\end{thm}
\begin{proof}
\proofstep{all:swap}{all:uniform}
Corollary~\ref{thm:prefix-pair} gives $\overline{r_v}=\overline{r_w}$
for every distinct $v,w$. For $s=1$ the assertion is automatic.

\proofstep{all:uniform}{all:degree}
Let $c\in R$ be the common residue. Then
$a=c\mathbf1^{\mathsf T}$ is the row of Lemma~\ref{lem:prefix-row}.
Its normalization gives $\overline{s}c=1$. In particular $c$ is a
unit of $R$, including when $d-1$ is composite. The $v$-entry of
$a\overline M=a$ is $c\,\overline{b(v)}=c$.
Multiplying by $c^{-1}$ gives $\overline{b(v)}=1$.

\proofstep{all:degree}{all:uniform}
The column-sum condition says $\mathbf1^{\mathsf T}\overline M
=\mathbf1^{\mathsf T}$. Iterating to level $L$ and using
\eqref{eq:folding-counts} yields
\[
 \mathbf1^{\mathsf T}
 =\mathbf1^{\mathsf T}\overline M^{\,L}
 =(\mathbf1^{\mathsf T}\mathbf1)\overline{\mathbf r}
 =\overline{s}\,\overline{\mathbf r}.
\]
Thus $\overline{s}\,\overline{r_v}=1$ for every $v$. Since $G^0$ is
nonempty, $\overline{s}$ has an inverse in $R$, so $\gcd(s,d-1)=1$
and $\overline{r_v}=\overline{s}^{\,-1}$ for all $v$.

\proofstep{all:uniform}{all:swap}
Apply Corollary~\ref{thm:prefix-pair} to each distinct pair.
\end{proof}

\begin{rem}
Theorem~\ref{thm:prefix-all} proves the equivalence of
\textup{(\ref{a:ind})} and \textup{(\ref{a:ex})} in Theorem~A.
By \eqref{eq:folding-counts}, $b(v)$ is the number of one-step
preimages of any path starting at $v$. Thus the condition on the
entire partition of $\fullshift$ can also be written as
\[
 |\sigma_G^{-1}(\{x\})|\equiv1\pmod{d-1}
 \qquad\text{for every }x\in X_G.
\]

For the example in Subsection~\ref{sub:blocks}, $d=3$ and
$(r_A,r_B,r_C)=(6,2,1)$. Hence $a=(0,0,1)$ over $\Z/2\Z$,
and Corollary~\ref{thm:prefix-pair} singles out $P_A,P_B$ as the
only pair admitting a transposition. The vertices $A,B$ have indegrees $6,1$,
respectively, with different residues modulo $2$. For the entire partition,
Theorem~\ref{thm:prefix-all} applies directly. The indegrees $(6,1,2)$ do not all have residue $1$,
so the three blocks cannot be interchanged pairwise.
\end{rem}

\section{The graph algebra and the induced maps on \texorpdfstring{$K_0$}{K0}}\label{sec:algebra}

Let $G$ be the graph obtained from a folding on $d\ge2$ letters,
and put $s=|G^0|$ and $R=\Z/(d-1)\Z$. The isomorphism and factor
maps in Subsections~\ref{sub:matrix-algebra} and~\ref{sub:fold-factors}
also apply for $d=2$. In Subsection~\ref{sub:fold-squares} we assume $d\ge3$. Fix the synchronizing level $L$ and its synchronizing map $t$ as in
Subsection~\ref{sub:foldings}. Theorem~\ref{thm:prefix-all} characterizes
pairwise transpositions of the sets $P_v$ by the indegrees of $G$.
Here we show that the same indegree condition characterizes commutativity
of the two $K_0$-diagrams for the factor maps of the folding.
The calculation uses the classes of the cylinder projections in the graph algebra.

\subsection{An isomorphism with a matrix algebra over \texorpdfstring{$\mathcal O_d$}{Od}}\label{sub:matrix-algebra}

We identify $\OG$ with $M_s(\Od)$ in a way that records both the
vertex projections and the canonical diagonal. This determines the
class of a single path cylinder and the class of the unit.

Index the standard matrix units $E_{ij}$ of $M_s(\C)$ by $i,j\in G^0$,
and put $\mathcal D_s=\operatorname{span}\{E_{ii}:i\in G^0\}$.
For $i\in G^0$ and a word $u$, write $\mu(i,u)$ for the unique path
with source $i$ and label $u$, with $\mu(i,\varepsilon)=i$. Let
$s_a$ ($a\in\mathcal A$) be the canonical Cuntz isometries of $\Od$ and put
$s_u=s_{u_1}\cdots s_{u_m}$ for a word $u=u_1\cdots u_m$, with
$s_\varepsilon=1$. Write $\gamma^G,\gamma^d$ for the gauge actions
on $\OG,\Od$, respectively. At $z\in\T$ these actions multiply every
edge or letter generator by $z$ and fix the vertex projections.
We identify $M_s(\Od)$ with $M_s(\C)\otimes\Od$.

\begin{prop}\label{prop:fold-matrix}
There is a unital $*$-isomorphism $\rho:\OG\to M_s(\Od)$ given by
\begin{equation}\label{eq:fold-matrix-generators}
 \rho(p_i)=E_{ii}\otimes1,\qquad
 \rho(S_e)=E_{s(e),r(e)}\otimes s_{\ell(e)} .
\end{equation}
It intertwines $\gamma^G$ with $\id\otimes\gamma^d$ and satisfies
\begin{equation}\label{eq:fold-matrix-diagonal}
 \rho\big(j_{X_G}(C(X_G))\big)
   =\mathcal D_s\otimes j_{\fullshift}(C(\fullshift)).
\end{equation}
\end{prop}
\begin{proof}
Write $Q_i=E_{ii}\otimes1$ and $T_e=E_{s(e),r(e)}\otimes s_{\ell(e)}$.
Then $T_e^*T_e=Q_{r(e)}$. For distinct edges with the same source, the labels
are distinct and the range projections are orthogonal. For different
sources the range projections lie in orthogonal matrix corners. Since
each vertex emits one edge of each label,
\[
 \sum_{s(e)=i}T_eT_e^*
   =E_{ii}\otimes\sum_{a\in\mathcal A}s_as_a^*=Q_i.
\]
The universal property for \eqref{eq:graph-relations} gives a
$*$-homomorphism $\rho$, unital since $\sum_i Q_i=1$. The formulas on
generators give gauge equivariance. Each $Q_i$ is nonzero, so $\rho$
is injective by the gauge-invariant uniqueness theorem
\cite[Theorem~2.2]{raeburn}.

To prove surjectivity, multiply the images of the successive edges of
$\mu(i,u)$. The adjacent matrix indices agree, so
\eqref{eq:fold-matrix-generators} gives
\begin{equation}\label{eq:fold-matrix-path}
 \rho(S_{\mu(i,u)})=E_{i,\delta(i,u)}\otimes s_u .
\end{equation}
For $|u|=L$, Lemma~\ref{lem:folding-form} gives
$\delta(i,u)=t(u)$ for every $i$.
Thus the two paths read from $i$ and $j$ have the same terminal vertex,
and
\[
 \sum_{u\in\mathcal A^L}
   \rho\big(S_{\mu(i,u)}S_{\mu(j,u)}^*\big)
   =E_{ij}\otimes\sum_{u\in\mathcal A^L}s_us_u^*
   =E_{ij}\otimes1 .
\]
The last equality is \eqref{eq:path-refinement} iterated $L$ times
in the one-vertex graph $E_d$.
For $a\in\mathcal A$, multiplying $\rho(S_{\mu(i,a)})$ on the right by
$E_{\delta(i,a),j}\otimes1$ now gives $E_{ij}\otimes s_a$. These
elements and the matrix units generate $M_s(\Od)$, so $\rho$ is onto.

Finally, \eqref{eq:fold-matrix-path} sends the range projection of
$S_{\mu(i,u)}$ to $E_{ii}\otimes s_us_u^*$. Every pair $(i,u)$ occurs.
Lemma~\ref{lem:diagonal-embedding} identifies the closed spans of these
projections on the two sides, giving \eqref{eq:fold-matrix-diagonal}.
\end{proof}

\begin{cor}\label{cor:fold-scaled-k}
All vertex projections have the same class $e\in K_0(\OG)$, and
\begin{equation}\label{eq:fold-scaled-k}
 K_0(\OG)=Re,\qquad [1_{\OG}]=se,\qquad K_1(\OG)=0 ,
\end{equation}
where $e$ has order $d-1$. For every finite path $\mu$, including a
vertex path, $c_{X_G}([\mu])=e$. For every word $u$,
$c_{\fullshift}([u])=\bar1$, where $\bar1=[1_{\Od}]$.
\end{cor}
\begin{proof}
By \cite[Theorems~3.7 and~3.8]{cuntz}, $K_0(\Od)$ is cyclic of
order $d-1$ with generator $\bar1$, and $K_1(\Od)=0$.
For each $i$, the corner embedding $\iota_i(x)=E_{ii}\otimes x$
induces isomorphisms on $K_0$ and $K_1$ by matrix stability
\cite[Propositions~4.3.8 and~8.2.8]{rll}. By
\eqref{eq:fold-matrix-generators},
\[
 \rho_{*0}([p_i])=[E_{ii}\otimes1]=(\iota_i)_{*0}(\bar1).
\]
These projections are Murray--von Neumann equivalent via
$E_{ij}\otimes1$. Their common class therefore gives a generator $e$
in $K_0(\OG)$ under $\rho_{*0}^{-1}$, and $[1_{\OG}]=\sum_i[p_i]=se$.
The isomorphisms in degree one give $K_1(\OG)=0$. The cylinder formulas follow from
\eqref{eq:path-class} for $G$ and for the one-vertex graph $E_d$.
\end{proof}

\subsection{Factor maps and cylinder preimages}\label{sub:fold-factors}

To compute the maps on cylinder classes, we describe the inverse image
of a cylinder under each of the two factor maps. Use the homeomorphism $\Theta:X_G\to G^0\times\fullshift$ of
Lemma~\ref{lem:path-coordinates}. Define $p:X_G\to\fullshift$ and
$q:\fullshift\to X_G$ by
\begin{equation}\label{eq:fold-pq}
 p(\Theta^{-1}(i,z))=z,\qquad
 q(uz)=\Theta^{-1}(t(u),z)
 \quad(i\in G^0,\ u\in\mathcal A^L,\ z\in\fullshift).
\end{equation}
Thus $p$ reads all the labels of a path. For $q$, the first $L$ letters
determine the initial vertex and the remaining letters determine the path.

\begin{lem}\label{lem:fold-factors}
Both $p$ and $q$ are factor maps, that is, continuous surjections
commuting with the respective shifts. They satisfy
\begin{equation}\label{eq:fold-compositions}
 p q=\sigma_d^L,\qquad q p=\sigma_G^L .
\end{equation}
For every word $w$ and every vertex $v$,
\begin{equation}\label{eq:fold-cylinder-preimages}
 \begin{aligned}
 p^{-1}([w])&=\bigsqcup_{i\in G^0}[\mu(i,w)],\\
 q^{-1}([\mu(v,w)])&=\bigsqcup_{\substack{u\in\mathcal A^L\\t(u)=v}}[uw].
 \end{aligned}
\end{equation}
\end{lem}
\begin{proof}
In the coordinates $\Theta$, $p$ is projection onto the second factor.
It is a continuous surjection and \eqref{eq:shift-coordinates} gives
$p\sigma_G=\sigma_dp$. On each clopen cylinder $[u]$, $|u|=L$, the
map $q$ is $y\mapsto\Theta^{-1}(t(u),\sigma_d^Ly)$, so $q$ is
continuous. Since $t$ is onto by Lemma~\ref{lem:folding-form},
for each $\Theta^{-1}(v,z)$ we may choose $u$ with $t(u)=v$.
Then $uz$ is a preimage under $q$.

By definition $pq(uz)=z$. On a path with coordinates $(i,uz)$, both
$qp$ and $\sigma_G^L$ have coordinates $(t(u),z)$ by
\eqref{eq:shift-coordinates}. This proves \eqref{eq:fold-compositions}.
Consequently
\[
 (q\sigma_d)p=qp\sigma_G=\sigma_G^{L+1}=(\sigma_G q)p.
\]
Surjectivity of $p$ gives $q\sigma_d=\sigma_Gq$.
By Lemma~\ref{lem:path-coordinates},
$\Theta([\mu(i,w)])=\{i\}\times[w]$. Thus the first union in
\eqref{eq:fold-cylinder-preimages} is the inverse image of $[w]$
under the coordinate projection. For the second, $q(uz)\in[\mu(v,w)]$
if and only if $t(u)=v$ and $z\in[w]$. Distinct length-$L$ words $u$
give disjoint cylinders $[uw]$. Both descriptions include the empty
word $w=\varepsilon$.
\end{proof}

The two unions have $s$ and $r_v$ summands, respectively, by
\eqref{eq:two-blocks}. By additivity and Corollary~\ref{cor:fold-scaled-k},
\eqref{eq:fold-cylinder-preimages} yields
\begin{equation}\label{eq:fold-coefficient-counts}
 c_{X_G}(p^{-1}([w]))=se,\qquad
 c_{\fullshift}(q^{-1}([\mu(v,w)]))=r_v\bar1 .
\end{equation}
The first class is independent of $w$. The second depends only on the
initial vertex $v$. In particular $q^{-1}(V_v)=P_v$.

\subsection{The indegree criterion for the induced maps on \texorpdfstring{$K_0$}{K0}}\label{sub:fold-squares}

Assume $d\ge3$. We now ask when a single isomorphism
$k:K_0(\Od)\to K_0(\OG)$ preserving the unit class makes the two
diagrams for $p$ and $q$ commute. Under the pointwise-rank identification
of Lemma~\ref{lem:clopen-classes}, a continuous map $h:Z\to W$
acts by $h^*:C(W,\Z)\to C(Z,\Z)$, $f\mapsto f\circ h$.
Thus the required identities are
\begin{equation}\label{eq:fold-squares}
 \begin{aligned}
 k\circ(j_{\fullshift})_{*0}&=(j_{X_G})_{*0}\circ p^*
       &&\text{on } C(\fullshift,\Z),\\
 k^{-1}\circ(j_{X_G})_{*0}&=(j_{\fullshift})_{*0}\circ q^*
       &&\text{on } C(X_G,\Z).
 \end{aligned}
\end{equation}
Both shift algebras have trivial $K_1$ by Corollary~\ref{cor:fold-scaled-k}.
On an indicator function, $h^*1_U=1_{h^{-1}(U)}$, so
\eqref{eq:fold-coefficient-counts} computes the right-hand sides on cylinders.

\begin{prop}\label{prop:fold-square-criterion}
The following conditions are equivalent.
\begin{enumerate}[\rm(i)]
\item\label{sq:maps} There exists a unit-preserving isomorphism $k$
satisfying \eqref{eq:fold-squares}.
\item\label{sq:arithmetic} $s r_v\equiv1\pmod{d-1}$ for every $v\in G^0$.
\item\label{sq:degree} $b(v)\equiv1\pmod{d-1}$ for every $v\in G^0$.
\end{enumerate}
Whenever they hold, $k$ is uniquely determined by $k(\bar1)=se$.
\end{prop}
\begin{proof}
\proofstep{sq:maps}{sq:arithmetic}
Unit preservation gives $k(\bar1)=se$. Apply the second identity in
\eqref{eq:fold-squares} to $1_{V_v}$. Equations~\eqref{eq:path-class}
and~\eqref{eq:fold-coefficient-counts}, with $w=\varepsilon$, give
$k^{-1}(e)=r_v\bar1$.
Applying $k$ gives $e=s r_v e$. Since $e$ has order $d-1$, condition
\textup{(\ref{sq:arithmetic})} follows.

\proofstep{sq:arithmetic}{sq:maps}
For any vertex $v$, $\overline{s}\,\overline{r_v}=1$ in $R$. Thus
$\overline{s}$ is invertible and $k(m\bar1)=ms e$ defines a
unit-preserving isomorphism. On a word cylinder $[w]$ the first identity
in \eqref{eq:fold-squares} reads $k(\bar1)=se$. On a path cylinder
$[\mu(v,w)]$ the second reads $k^{-1}(e)=r_v\bar1$, which holds
since $k(r_v\bar1)=s r_v e=e$. In each case the asserted classes
are those in \eqref{eq:fold-coefficient-counts}.
As recalled in Subsection~\ref{sub:compatibility}, the cylinder
indicators generate both groups of continuous integer-valued functions. Hence both identities extend by
additivity to the indicated groups. Since $\bar1$ generates $K_0(\Od)$,
the value on $\bar1$ determines $k$.

\proofstep{sq:arithmetic}{sq:degree}
All $\overline{r_v}$ equal $\overline{s}^{\,-1}$. Apply
\textup{(\ref{all:uniform}) $\Rightarrow$ (\ref{all:degree})} of
Theorem~\ref{thm:prefix-all}.

\proofstep{sq:degree}{sq:arithmetic}
Theorem~\ref{thm:prefix-all} gives $\overline{r_v}=\overline{s}^{\,-1}$
for every $v$.
\end{proof}

Since $se$ generates $Re$ precisely when $\gcd(s,d-1)=1$, this
coprimality condition characterizes the existence of a unit-preserving
isomorphism $k$. For any such $k$, the forward identity holds on every
word cylinder by \eqref{eq:fold-coefficient-counts}. In the reverse
identity the left-hand side has value $k^{-1}(e)$ and the right-hand
side has value $r_v\bar1$ on every $[\mu(v,w)]$, including $V_v$
when $w=\varepsilon$. It therefore suffices to check the reverse
identity on the vertex cylinders. By Theorem~\ref{thm:prefix-all},
such an isomorphism makes both diagrams commute precisely when the sets $P_v$ admit
pairwise transpositions in $V_{d,1}$.

The maps $p,q$ in \eqref{eq:fold-squares} are factor maps.
Definition~\ref{def:ak} requires sequences of homeomorphisms. Section~\ref{sec:real} uses the
inverse images in \eqref{eq:fold-cylinder-preimages} to realize the
action of each factor map on successive finite cylinder partitions.

\section{Approximate \texorpdfstring{$K$}{K}-conjugacy with the full shift}\label{sec:real}

By Proposition~\ref{prop:fold-square-criterion}, the indegree condition
determines when the factor maps $p$ and $q$ satisfy the two coefficient
identities. To prove Theorem~A we must realize these identities and the
shift approximations along the same homeomorphisms. Conversely, given
homeomorphisms satisfying Definition~\ref{def:ak}, we
recover the residue of each indegree modulo $d-1$ by comparing a vertex
cylinder with its one-step inverse image.

\subsection{Homeomorphisms and cylinder partitions}\label{sub:real-partitions}

At each word length we prescribe the inverse images of all cylinders.
Their overlaps control the shift approximation, and their projection
classes determine whether a single unitary realizes the prescribed
action on the finite-dimensional algebra of cylinder functions.

In this subsection let $Z$ and $W$ be one-sided subshifts whose spaces
are Cantor spaces, with the canonical embeddings $j_Z,j_W$ of
Subsection~\ref{sub:compatibility}. Write $\mathcal L_n(W)$ for the
words of length $n$ occurring in $W$.
The associated nonempty cylinders form the partition $\mathcal P_n(W)$.
For a compatible metric $d_W$, write the mesh of this partition as
\[
 \delta_n(W)=\max_{w\in\mathcal L_n(W)}\operatorname{diam}_{d_W}[w].
\]
Then $\delta_n(W)\to0$ since the cylinder partitions refine one another
and form a basis on the compact space $W$.

\begin{lem}\label{lem:real-partition}
Fix $n\ge2$ and a clopen partition $(D_w)_{w\in\mathcal L_n(W)}$ of
$Z$ into nonempty sets. Suppose that whenever $D_w\cap\sigma_Z^{-1}D_{w'}$
is nonempty, the last $n-1$ symbols of $w$ agree with the first $n-1$
symbols of $w'$. Then there is a homeomorphism $h:Z\to W$ such that
\begin{equation}\label{eq:real-overlap}
 h^{-1}([w])=D_w\quad(w\in\mathcal L_n(W)),\qquad
 \sup_{z\in Z}d_W\big(h(\sigma_Zz),\sigma_W(hz)\big)\le\delta_{n-1}(W).
\end{equation}
Suppose in addition that $\theta:\mathcal O_W\to\mathcal O_Z$ is a
unital $*$-isomorphism, $\mathcal O_Z$ is simple and purely infinite, and
\begin{equation}\label{eq:real-atom-class}
 c_Z(D_w)=\theta_{*0}(c_W([w]))\qquad(w\in\mathcal L_n(W)).
\end{equation}
For every homeomorphism $h$ with these cylinder preimages there is a
unitary $u\in\mathcal O_Z$ satisfying
\begin{equation}\label{eq:real-unitary}
 (\Ad u\circ\theta)(j_W(f))=j_Z(f\circ h)
 \qquad\big(f\in\operatorname{span}\{1_{[w]}:w\in\mathcal L_n(W)\}\big).
\end{equation}
\end{lem}
\begin{proof}
Each $D_w$ and each $[w]$ is a nonempty clopen subset of a Cantor
space, hence is itself homeomorphic to the Cantor space by Brouwer's
theorem \cite[Theorem~7.4]{kechris}. Choose a homeomorphism from $D_w$
onto $[w]$ for each $w$. Their union is a homeomorphism $h:Z\to W$
because both partitions are finite and clopen. For any such $h$,
if $z\in D_w$ and
$\sigma_Z z\in D_{w'}$, then $\sigma_W(hz)$ starts with the last $n-1$
symbols of $w$, and $h(\sigma_Z z)$ starts with the first $n-1$
symbols of $w'$. By hypothesis these two points belong to the same
length-$(n-1)$ cylinder, giving \eqref{eq:real-overlap}.

Put $P_w=\theta(j_W(1_{[w]}))$ and $Q_w=j_Z(1_{D_w})$. Each family
consists of mutually orthogonal projections summing to the unit.
They are nonzero since the atoms are nonempty and the canonical
embeddings and $\theta$ are injective.
In the unital simple purely infinite algebra $\mathcal O_Z$, nonzero
projections of equal $K_0$-class are Murray--von Neumann equivalent
\cite[Exercise~5.7(ii)]{rll}.
Thus \eqref{eq:real-atom-class} gives $v_w\in\mathcal O_Z$ with
$v_w^*v_w=P_w$ and $v_wv_w^*=Q_w$. Since $v_w=Q_wv_wP_w$,
orthogonality gives
\[
 u=\sum_w v_w,\qquad u^*u=\sum_wP_w=1,\qquad
 uu^*=\sum_wQ_w=1,\qquad uP_wu^*=Q_w .
\]
This proves \eqref{eq:real-unitary} on the indicators and hence on
their linear span.
\end{proof}

\begin{cor}\label{cor:real-factor}
Let $\theta:\mathcal O_W\to\mathcal O_Z$ be a unital $*$-isomorphism,
where $\mathcal O_Z$ is simple and purely infinite. For each $n\ge1$
suppose there is a continuous map $F_n:Z\to W$ commuting with the
shifts such that, for every $w\in\mathcal L_{n+1}(W)$,
\[
 F_n^{-1}([w])\ne\varnothing,\qquad
 c_Z(F_n^{-1}([w]))=\theta_{*0}(c_W([w])).
\]
Then there are homeomorphisms $h_n:Z\to W$ and unital
$*$-isomorphisms $\theta_n:\mathcal O_W\to\mathcal O_Z$ such that
$[\theta_n]_{KL}=[\theta]_{KL}$ and
\begin{align}
 \sup_{z\in Z}d_W\big(h_n(\sigma_Z z),\sigma_W(h_nz)\big)&\longrightarrow0,
       \label{eq:real-uniform}\\
 \|\theta_n(j_W(f))-j_Z(f\circ h_n)\|&\longrightarrow0
       \qquad(f\in C(W)).\label{eq:real-norm}
\end{align}
\end{cor}
\begin{proof}
By \cite[Lemma~3.2 and its proof]{shimomura}, applied to $F_n$ and
$\mathcal P_{n+1}(W)$, there are homeomorphisms $h_n:Z\to W$ with
\[
 h_n^{-1}([w])=F_n^{-1}([w])\quad(w\in\mathcal L_{n+1}(W)),\qquad
 h_n\sigma_Z h_n^{-1}\longrightarrow\sigma_W\quad\text{uniformly}.
\]
Indeed, each image $F_n(Z)$ meets every atom at that level, as required
by Shimomura's lemma. The uniform limit gives \eqref{eq:real-uniform}.
The class assumption yields \eqref{eq:real-atom-class} for these
preimages, and $F_n\sigma_Z=\sigma_W F_n$ gives the required overlap.
Use the unitary assertion of Lemma~\ref{lem:real-partition} for $h_n$ and
put $\theta_n=\Ad u_n\circ\theta$. Unitarily equivalent homomorphisms
have the same $KK$-class \cite[Sections~17.1--17.2]{blackadar}, hence
$[\theta_n]_{KL}=[\theta]_{KL}$.

Fix $f\in C(W)$ and $\varepsilon>0$. By the density recalled in
Subsection~\ref{sub:compatibility}, choose $g$ constant on cylinders of
some length $m$ with $\|f-g\|<\varepsilon$. For $n+1\ge m$,
\eqref{eq:real-unitary} applies to $g$ by refinement. Contractivity gives
$\|\theta_n(j_W(f))-j_Z(f\circ h_n)\|\le2\|f-g\|<2\varepsilon$.
This proves \eqref{eq:real-norm} for the same $h_n$.
\end{proof}

\subsection{Classes of one-step inverse images}\label{sub:real-inverse}

For the converse we use only the homeomorphisms from $\fullshift$ to
$X_G$. The relevant property of the full shift is the identity
\begin{equation}\label{eq:real-full-invariance}
 c_{\fullshift}(\sigma_d^{-1}U)=c_{\fullshift}(U)
 \qquad(U\subseteq\fullshift\text{ clopen}).
\end{equation}
Indeed, $\sigma_d^{-1}[w]=\bigsqcup_{a\in\mathcal A}[aw]$ and each
cylinder has class $\bar1$ by Corollary~\ref{cor:fold-scaled-k}.
Thus the inverse image has class $d\bar1=\bar1$. A disjoint cylinder
decomposition and Lemma~\ref{lem:clopen-classes} give
\eqref{eq:real-full-invariance} for every clopen $U$.

For $v\in G^0$, put
\begin{equation}\label{eq:real-defect}
 \Delta_v=c_{X_G}(\sigma_G^{-1}V_v)-c_{X_G}(V_v)=(b(v)-1)e .
\end{equation}
Here $e$ is the generator of Corollary~\ref{cor:fold-scaled-k}.
Indeed $\sigma_G^{-1}V_v=\bigsqcup_{r(\alpha)=v}[\alpha]$, where
$\alpha$ ranges over the incoming edges. Each of the $b(v)$ edge
cylinders and the vertex cylinder $V_v$ have class $e$.

\begin{lem}\label{lem:real-transport}
Let $X$ be a one-sided subshift and let $\psi_n:\fullshift\to X$
be homeomorphisms such that
$\sup_{y\in\fullshift}d_X(\psi_n(\sigma_d y),\sigma_X(\psi_n y))\to0$.
Suppose $k:K_0(\Od)\to K_0(\mathcal O_X)$ is a group isomorphism
and for each clopen $B\subseteq X$,
$c_{\fullshift}(\psi_n^{-1}B)=k^{-1}(c_X(B))$ for all sufficiently large $n$.
Then $c_X(\sigma_X^{-1}B)=c_X(B)$ for every clopen $B\subseteq X$.
In particular, if $X=X_G$ for a folded automaton $G$, then $\Delta_v=0$ for all $v$.
\end{lem}
\begin{proof}
Fix $B$. The assertion is immediate for $B=\varnothing,X$.
Otherwise $B$ and its complement are disjoint compact sets at
positive distance. By the assumed uniform limit, for all large $n$
the points $\psi_n(\sigma_d y)$ and $\sigma_X(\psi_n y)$ have the
same $B$-membership for every $y\in\fullshift$. Hence
\[
 \psi_n^{-1}(\sigma_X^{-1}B)=\sigma_d^{-1}(\psi_n^{-1}B).
\]
The set $\sigma_X^{-1}B$ is clopen by continuity. Take one $n$ for
which the displayed set equality and the two class identities, for
$B$ and for $\sigma_X^{-1}B$, all hold. This is possible since each
holds for all sufficiently large $n$. Using \eqref{eq:real-full-invariance}
on the clopen set $\psi_n^{-1}B$, we obtain
\[
 \begin{aligned}
 k^{-1}(c_{X}(\sigma_X^{-1}B))
 &=c_{\fullshift}(\psi_n^{-1}(\sigma_X^{-1}B))
  =c_{\fullshift}(\sigma_d^{-1}(\psi_n^{-1}B))\\
 &=c_{\fullshift}(\psi_n^{-1}B)=k^{-1}(c_{X}(B)).
 \end{aligned}
\]
Injectivity of $k^{-1}$ gives the assertion. Take $B=V_v$ for the
last claim.
\end{proof}

\subsection{Proof of Theorem A}\label{sub:real-main}

For the ordered pair $(X_G,\fullshift)$, the sufficient condition
uses the specific factor maps of Subsection~\ref{sub:fold-factors}.
For the converse, Lemma~\ref{lem:real-transport} applies to any reverse
sequence supplied by Definition~\ref{def:ak}.

\begin{proof}[Proof of Theorem~A]
\proofstep{a:ind}{a:ak}
By Proposition~\ref{prop:fold-square-criterion}, there is a
unit-preserving isomorphism $k:K_0(\Od)\to K_0(\OG)$ satisfying
\eqref{eq:fold-squares} in both directions. Both $K_1$-groups vanish by
Corollary~\ref{cor:fold-scaled-k}, so $k$ extends uniquely to a graded
isomorphism. The algebras are unital, separable, nuclear, simple and
purely infinite by Proposition~\ref{prop:fold-matrix} and the properties
of $\Od$ \cite{cuntz1977,cuntz}. They satisfy the universal coefficient
theorem as recalled in Subsection~\ref{sub:compatibility}. The unital
Kirchberg--Phillips theorem \cite[Theorem~4.2.4]{phillips} gives a
$*$-isomorphism $\Psi:\Od\to\OG$ with $\Psi_{*0}=k$. Put
$\Phi=\Psi^{-1}$ and $\kappa=[\Psi]_{KL}$.

Apply the first identity of \eqref{eq:fold-squares} to each cylinder
indicator on $\fullshift$. It gives the hypothesis of
Corollary~\ref{cor:real-factor} for the constant sequence $F_n=p$ and $\theta=\Psi$. The resulting
$\phi_n:X_G\to\fullshift$ and $\Psi_n:\Od\to\OG$ satisfy the first
limits in \eqref{eq:ak-uniform} and \eqref{eq:ak-norm}, with
$[\Psi_n]_{KL}=\kappa$. Apply the second identity to each cylinder
indicator on $X_G$ and use the same corollary with $F_n=q$ and $\theta=\Phi$. Both maps are onto by Lemma~\ref{lem:fold-factors}.
This gives $\psi_n:\fullshift\to X_G$ and $\Phi_n:\OG\to\Od$
satisfying the second limits, with
$[\Phi_n]_{KL}=[\Phi]_{KL}=\kappa^{-1}$. Thus the two sequences of homeomorphisms and the corresponding
algebra isomorphisms satisfy Definition~\ref{def:ak} with class $\kappa$.

\proofstep{a:ak}{a:ind}
Take $\psi_n$ from Definition~\ref{def:ak} and $k=\kappa_{*0}$.
By Corollary~\ref{cor:ak-finite}, \eqref{eq:finite-reverse} holds
along these same $\psi_n$. For each clopen $B\subseteq X_G$, take
$J=\langle 1_B,1_{\sigma_G^{-1}B}\rangle$ and $H=0$ in
Definition~\ref{def:finitecoeff}. Both class identities then hold for
every $n\ge N(0,J)$. This gives the hypothesis on clopen classes in
Lemma~\ref{lem:real-transport}, while \eqref{eq:ak-uniform} supplies
its dynamical hypothesis. Thus $\Delta_v=0$ for each $v$.
By \eqref{eq:real-defect} and the fact that $e$ has order $d-1$,
$b(v)\equiv1\pmod{d-1}$ for every vertex.

\proofstep{a:ind}{a:ex}
Apply implication \textup{(\ref{all:degree}) $\Rightarrow$
(\ref{all:swap})} of Theorem~\ref{thm:prefix-all}, via
\textup{(\ref{all:uniform})}.

\proofstep{a:ex}{a:ind}
Apply \textup{(\ref{all:swap}) $\Rightarrow$ (\ref{all:uniform})
$\Rightarrow$ (\ref{all:degree})} of the same theorem.
\end{proof}

At a fixed word length the shift approximation can be expressed by
the overlap of consecutive words. The next corollary characterizes
approximate $K$-conjugacy by partitions with this overlap and the
prescribed $K_0$-classes. The same $K_0$-isomorphism is used at every
length, while the partitions may be chosen separately at each length.

\begin{cor}\label{cor:real-equivalents}
For the ordered pair $(X_G,\fullshift)$, the following are equivalent.
\begin{enumerate}[\rm(i)]
\item\label{real:ak} $X_G$ and $\fullshift$ are approximately $K$-conjugate.
\item\label{real:both} There exist a unit-preserving graded isomorphism $\eta$ and
sequences of homeomorphisms $(\phi_n),(\psi_n)$ as in
Definition~\ref{def:finitecoeff} for which both the finite coefficient
$K$-condition and the two limits in \eqref{eq:ak-uniform} hold.
\item\label{real:delta} $\Delta_v=0$ for every $v\in G^0$.
\item\label{real:factors} Some unit-preserving isomorphism $k:K_0(\Od)\to K_0(\OG)$
satisfies \eqref{eq:fold-squares} for $p,q$.
\item\label{real:parts} There is a unit-preserving isomorphism $k:K_0(\Od)\to K_0(\OG)$
such that for each $n\ge2$ and each of the triples
\[
 (Z,W,\lambda)=(X_G,\fullshift,k),\qquad
 (Z,W,\lambda)=(\fullshift,X_G,k^{-1}),
\]
there is a clopen partition
$(D_w)_{w\in\mathcal L_n(W)}$ of $Z$ into nonempty sets satisfying the
overlap condition of Lemma~\ref{lem:real-partition} and
$c_Z(D_w)=\lambda(c_W([w]))$ for every $w$.
\end{enumerate}
\end{cor}
\begin{proof}
\proofstep{real:ak}{real:both}
Use Corollary~\ref{cor:ak-finite} with the homeomorphisms of
Definition~\ref{def:ak}.

\proofstep{real:both}{real:delta}
Apply Lemma~\ref{lem:real-transport} to the reverse homeomorphisms,
using \eqref{eq:finite-reverse} on the subgroup generated by $1_B$
for each clopen $B\subseteq X_G$.

\proofstep{real:delta}{real:factors}
By \eqref{eq:real-defect} and the order of $e$, the indegree
condition holds. Apply Proposition~\ref{prop:fold-square-criterion}.

\proofstep{real:factors}{real:parts}
Use $D_w=p^{-1}([w])$ and $D_w=q^{-1}([w])$ in the respective directions.
The class equalities are \eqref{eq:fold-squares} on the indicators.
Each family has the required overlap since $p,q$ commute with the
shifts, as in Corollary~\ref{cor:real-factor}.

\proofstep{real:parts}{real:both}
Extend $k$ by the unique map on the zero $K_1$-groups
(Corollary~\ref{cor:fold-scaled-k}) to obtain $\eta$.
For each $n\ge1$ use the first assertion of Lemma~\ref{lem:real-partition}
in both directions at length $n+1$. It gives $\phi_n,\psi_n$ with the
prescribed cylinder preimages and both limits in \eqref{eq:ak-uniform}.
To verify \eqref{eq:atom-diagrams}, fix finite clopen partitions on the
two shifts. Choose $m$ so that both length-$m$ cylinder partitions
refine the respective clopen partitions, as in
Lemma~\ref{lem:clopen-classes}. For each $n+1\ge m$, every fixed atom is
a disjoint union of length-$(n+1)$ cylinders. Adding the prescribed
class identities over that union gives \eqref{eq:atom-diagrams} for the
same $\phi_n,\psi_n$. Lemma~\ref{lem:atoms} now gives the finite
coefficient $K$-condition.

\proofstep{real:factors}{real:ak}
Proposition~\ref{prop:fold-square-criterion} gives
\textup{(\ref{a:ind})} of Theorem~A. The implication
\textup{(\ref{a:ind}) $\Rightarrow$ (\ref{a:ak})} proved above gives
the conclusion.
\end{proof}

The homeomorphisms in \textup{(\ref{real:both})} carry both the
coefficient identities and the shift approximations. Section~\ref{sec:examples}
considers the finite coefficient $K$-condition of Definition~\ref{def:finitecoeff}
on its own.

\section{The finite coefficient \texorpdfstring{$K$}{K}-condition and examples}\label{sec:examples}

Theorem~A requires the homeomorphisms in the coefficient diagrams to
approximate the shifts. We now determine when the diagrams alone can
be realized. For a fixed alphabet size $d$, the criterion compares the greatest common divisors of the vertex counts with $d-1$. We then compare it with the indegree criterion, examine the
additional information in the gauge fixed-point algebra, and treat a
two-vertex shift outside the class of foldings.

\subsection{The finite coefficient \texorpdfstring{$K$}{K}-condition}\label{sub:six-static}

We first compute the determinant needed for continuous orbit equivalence.
Together with the unit class computed in Corollary~\ref{cor:fold-scaled-k},
it gives a criterion for two foldings on the same alphabet.

Write $M_G$ for the vertex adjacency matrix of $G$ in the source-row
convention. To apply the orbit equivalence theorem to graphs with
parallel edges, we also use the edge matrix
$B_G=(1_{\{r(e)=s(f)\}})_{e,f\in G^1}$. Its one-sided topological
Markov shift is $X_G$. The standard identification
$\mathcal O_{B_G}\cong C^*(G)$ sends the generator indexed by $e$
to $S_e$ \cite[Remark~2.8]{raeburn}, using our path convention.
Thus the diagonal projections of words map to the corresponding path
projections. The same identification applies to any finite graph with
no sinks or sources used below.

\begin{lem}\label{lem:six-det}
For the graph $G$ obtained from a folding on $d\ge2$ letters, with $s=|G^0|$,
\[
 \det(\lambda I_s-M_G)=\lambda^{s-1}(\lambda-d),\qquad
 \det(I-B_G)=\det(I-M_G)=1-d .
\]
The matrix $B_G$ is irreducible and is not a permutation matrix.
\end{lem}
\begin{proof}
By \eqref{eq:folding-counts}, $M_G\mathbf1=d\mathbf1$ and the image
of $M_G^L=\mathbf1\mathbf r$ is the line spanned by $\mathbf1$ over
$\mathbb Q$. On the quotient by that line the induced matrix is
nilpotent. This gives the characteristic polynomial, also for $s=1$.
Set $S_{v,e}=1_{\{s(e)=v\}}$ and $T_{e,v}=1_{\{r(e)=v\}}$.
Then $M_G=ST$ and $B_G=TS$, so Sylvester's determinant identity gives
$\det(I-B_G)=\det(I-M_G)$. Each entry of $M_G^L$ is positive since
the synchronizing map is onto. Thus a path joins the range of any
edge to the source of any other, and $B_G$ is irreducible. Each row of
$B_G$ contains exactly $d$ ones.
\end{proof}

For shifts $(X,\sigma_X)$ and $(Y,\sigma_Y)$, continuous orbit
equivalence means that there is a homeomorphism $h:X\to Y$ and
continuous functions $k_X,l_X:X\to\mathbb N_0$ and
$k_Y,l_Y:Y\to\mathbb N_0$ such that
\begin{equation}\label{eq:six-orbit}
 \begin{aligned}
 \sigma_Y^{k_X(x)}h(\sigma_X x)&=\sigma_Y^{l_X(x)}h(x) &&(x\in X),\\
 \sigma_X^{k_Y(y)}h^{-1}(\sigma_Y y)&=\sigma_X^{l_Y(y)}h^{-1}(y) &&(y\in Y).
 \end{aligned}
\end{equation}
We use the one-sided convention of \cite{matmat}. The following statement includes Theorem~B by
taking one graph to be the one-vertex presentation of $\fullshift$.

\begin{thm}\label{thm:six-static}
Let $G,H$ be the quotient graphs of two foldings on the same alphabet of $d\ge3$ letters, with
$s=|G^0|$ and $t=|H^0|$. The following conditions are equivalent.
\begin{enumerate}[\rm(i)]
\item\label{six:fc} $X_G$ and $X_H$ satisfy the finite coefficient $K$-condition.
\item\label{six:gcd} $\gcd(s,d-1)=\gcd(t,d-1)$.
\item\label{six:iso} There is a unital $*$-isomorphism $\mathcal O_G\cong\mathcal O_H$.
\item\label{six:coe} $X_G$ and $X_H$ are continuously orbit equivalent.
\end{enumerate}
When these conditions hold, the two coefficient identities can be
realized on the entire groups of continuous integer-valued functions
by a single homeomorphism and its inverse.
\end{thm}
\begin{proof}
Write $m=d-1$. By Corollary~\ref{cor:fold-scaled-k} the respective
$K_0$-groups are cyclic of order $m$, with unit classes $s e_G$ and
$t e_H$ for generators $e_G,e_H$. Both $K_1$-groups vanish.

\proofstep{six:fc}{six:gcd}
The isomorphism in Definition~\ref{def:finitecoeff} preserves the unit
class. The orders of the two unit classes are $m/\gcd(s,m)$ and $m/\gcd(t,m)$,
respectively, proving the equality.

\proofstep{six:gcd}{six:iso}
By \cite[Lemma~1]{abramssmith} applied to $\Z/m\Z$ with generator
$\bar1$, the equality of the two greatest common divisors gives an
automorphism carrying $s\bar1$ to $t\bar1$. This identifies the
two $K_0$-groups with their distinguished unit classes.
By Proposition~\ref{prop:fold-matrix} the algebras are $M_s(\Od)$
and $M_t(\Od)$. The classification of matrix algebras over Cuntz
algebras by pointed $K_0$ now gives the required isomorphism
\cite[Theorem~5.12]{aalp}.

\proofstep{six:iso}{six:coe}
By Lemma~\ref{lem:six-det}, $B_G,B_H$ are irreducible non-permutation
$\{0,1\}$-matrices with $\det(I-B_G)=\det(I-B_H)$. By
\textup{(\ref{six:iso})} their Cuntz--Krieger algebras are isomorphic.
These are the hypotheses for \cite[Theorem~3.6]{matmat} to give
continuous orbit equivalence.

\proofstep{six:coe}{six:fc}
By the same theorem there is a diagonal-preserving isomorphism
$\Theta:\mathcal O_H\to\mathcal O_G$. Its restriction to the
diagonals induces a homeomorphism $h:X_G\to X_H$ satisfying
\[
 \Theta(j_{X_H}(f))=j_{X_G}(f\circ h)\qquad(f\in C(X_H)).
\]
The inverse satisfies the corresponding identity for $h^{-1}$.
Take $\eta=\Theta_*$, $\phi_n=h$ and $\psi_n=h^{-1}$ in
Definition~\ref{def:finitecoeff}. Applying $K_0$ gives both
coefficient identities for every $n$ on the whole groups.
\end{proof}

\begin{proof}[Proof of Theorem~B]
For the one-vertex graph $H$ with $d$ loops, $X_H=\fullshift$,
$\mathcal O_H=\Od$ and $\gcd(|H^0|,d-1)=1$. Each assertion follows
from the corresponding implication of Theorem~\ref{thm:six-static}.
\end{proof}

\subsection{Examples of foldings}\label{sub:six-fold-examples}

The first pair of examples separates the two criteria with the
alphabet and the number of vertices fixed. The second example
separates approximate $K$-conjugacy from eventual conjugacy.

\begin{exm}\label{ex:six-four}
Let $d=4$. Partition the letters into two classes and send each vertex
to the class of the letter read. This gives a folding at level one.
Consider class sizes $(1,3)$ and $(2,2)$, respectively. In each case
both rows of the adjacency matrix are the row of class sizes. Each of these edge shifts satisfies the finite coefficient $K$-condition with $\Sigma_4^+$ by Theorem~B since $\gcd(2,3)=1$. The indegrees are
$(2,6)$ in the first case and $(4,4)$ in the second. By Theorem~A
only the second shift is approximately $K$-conjugate to $\Sigma_4^+$.
In the first folding the two cylinder counts have different residues
modulo $3$. In the second choose an involution $\pi$ of the alphabet that interchanges the two pairs. The transposition is $ix\mapsto\pi(i)x$.
\end{exm}

For a folded automaton $G$ write $\mathcal F_G=\mathcal O_G^{\gamma^G}$ for the
gauge fixed-point algebra and $\mathcal F_d=\Od^{\gamma^d}$. The gauge
equivariance in Proposition~\ref{prop:fold-matrix} gives
$\mathcal F_G\cong M_s(\mathcal F_d)$ by taking fixed points entrywise.
The Cuntz core $\mathcal F_d$ is the UHF algebra of type $d^\infty$
\cite{cuntz1977}. In particular the unique normalized trace $\tau_G$ has
\begin{equation}\label{eq:six-core}
 (\tau_G)_*(K_0(\mathcal F_G))=\tfrac1s\Z[1/d],\qquad
 \tau_G(j_{X_G}(1_{[\mu]}))=\frac1{s d^{|\mu|}} .
\end{equation}
Indeed the matrix stages have sizes $s d^n$, so the trace ranges are
the union of $(s d^n)^{-1}\Z$ by continuity of $K_0$ \cite{blackadar}.
The cylinder projection corresponds to $E_{ii}\otimes s_us_u^*$ for
$i=s(\mu)$ and the label word $u$ of $\mu$, a rank-one projection at
stage $|\mu|$ by \eqref{eq:fold-matrix-path}.

\begin{exm}\label{ex:six-five}
Let $d=3$ and partition the length-two words into the five classes
\[
 \{01,11,21\},\quad\{02,12,22\},\quad\{00\},\quad\{10\},\quad\{20\}.
\]
Two words in one class have the same last letter, so appending a
letter produces the same length-two suffix. This is a folding.
The five cardinalities $r_v$, namely $(3,3,1,1,1)$, are all odd. Thus the blocks admit
pairwise transpositions by Theorem~\ref{thm:prefix-all}, and the edge
shift is approximately $K$-conjugate to $\Sigma_3^+$ by Theorem~A.
For instance, subdivide $[00]$ into the three cylinders $[00a]$, $a=0,1,2$, and exchange them
in order with $[01],[11],[21]$, fixing the complement.

Eventual conjugacy requires a homeomorphism satisfying \eqref{eq:six-orbit} with $l_X=k_X+1$ and $l_Y=k_Y+1$. It would induce a
diagonal-preserving isomorphism intertwining the gauge actions
\cite[Theorem~4.1]{carlsenrout}, hence a unital isomorphism of the
fixed-point algebras. This would preserve their unique normalized
traces. Formula~\eqref{eq:six-core} gives ranges
$\tfrac15\Z[1/3]$ and $\Z[1/3]$, respectively. The first contains
$1/5$ and the second does not. Hence the shifts are not eventually
conjugate.
\end{exm}

In fact the two-sided label map is a conjugacy for every folding.
For a bi-infinite sequence of labels $(a_j)_{j\in\Z}$, the source
of the edge at time $j$ must be $t(a_{j-L}\cdots a_{j-1})$ by
synchronization \cite[Proposition~3.4]{bleak}. Taking its unique outgoing edge labelled $a_j$
gives a path by \eqref{eq:suffix-transition}. This finite-window rule
is the continuous shift-commuting inverse to the label map. Thus the
distinctions in these examples concern the one-sided shifts.

\subsection{Obstructions for irreducible edge shifts}\label{sub:six-general}

For a general irreducible edge shift, invariance of clopen $K_0$-classes
under inverse images gives an obstruction at each vertex. The examples
below compare this obstruction with the existence of fixed points and
with invariance of a trace on the diagonal.

Let $E$ be a finite strongly connected graph with no sinks.
Write $M$ for its source-row matrix and $b(v)=\sum_i M_{iv}$.
The graph $K$-theory formula \cite[Theorem~7.16]{raeburn}, translated
to our convention, gives
\begin{equation}\label{eq:six-graph-k}
 K_0(C^*(E))=\operatorname{coker}(I-M^{\mathsf T}),\qquad
 K_1(C^*(E))=\ker(I-M^{\mathsf T}).
\end{equation}
The standard generators correspond to $[p_v]$, with unit $\sum_v[p_v]$.
In particular, vertex relations give a presentation of $K_0$.

\begin{prop}\label{prop:six-ann}
The following conditions are equivalent.
\begin{enumerate}[\rm(i)]
\item\label{six:inv} $c_{X_E}(\sigma_E^{-1}U)=c_{X_E}(U)$ for every clopen $U$.
\item\label{six:ann} $(b(v)-1)K_0(C^*(E))=0$ for every vertex $v$.
\end{enumerate}
Approximate $K$-conjugacy with $\fullshift$ implies both conditions and
$b(v)\equiv1\pmod{d-1}$ at every vertex.
\end{prop}
\begin{proof}
For a finite path $\mu$ from $v$ to $w$, prepending any of the $b(v)$
edges ending at $v$ gives a cylinder of class $[p_w]$. Indeed, for
every finite path $\nu$ the graph partial isometry $S_\nu$ satisfies
$S_\nu S_\nu^*=j_{X_E}(1_{[\nu]})$ and $S_\nu^*S_\nu=p_{r(\nu)}$
\cite{raeburn}, as in \eqref{eq:path-class}. Consequently
\[
 c_{X_E}(\sigma_E^{-1}[\mu])-c_{X_E}([\mu])=(b(v)-1)[p_w].
\]
\proofstep{six:inv}{six:ann}
Fix $v$. Strong connectivity allows every $w$ to occur, and the
$[p_w]$ generate by \eqref{eq:six-graph-k}.
\proofstep{six:ann}{six:inv}
The displayed difference vanishes for every cylinder. Add over a
finite disjoint cylinder decomposition of $U$.

For the last assertion, apply Lemma~\ref{lem:real-transport} with $X=X_E$ and $k=\kappa_{*0}$. Its uniform limit is the second limit in \eqref{eq:ak-uniform}.
Its class identity is \eqref{eq:finite-reverse} on the subgroup
generated by $1_B$ for each clopen $B\subseteq X_E$, by Corollary~\ref{cor:ak-finite}. The resulting $K_0$-group
is cyclic of order $d-1$, so \textup{(\ref{six:ann})} gives the congruence.
\end{proof}

\begin{exm}\label{ex:six-loops}
For $j=1,2$ let $E_j$ be the directed graph with adjacency matrix
\[
 A_j=\begin{pmatrix}0&1&0\\j&0&1\\1&0&0\end{pmatrix}.
\]
Each $E_j$ is strongly connected and has cycles of lengths two and three. Thus its edge shift $X_{E_j}$ is mixing \cite[Chapter~4]{lindmarcus}. Write $e_i=[p_i]$. The
relations $e_1=e_2$, $e_2=je_1+e_3$, $e_3=e_1$ give $K_0=\Z/j\Z$
with unit $3e_1$. Since $\det(I-A_j)=-j$, $K_1=0$. The unit generates
for both values of $j$. The corresponding full $(j+1)$-shift has the
same pointed $K_0$ and determinant. Passing to edge matrices as in
Lemma~\ref{lem:six-det} preserves the determinant. The edge matrices
are irreducible and are not permutation matrices since both graphs
are strongly connected and have a branching vertex. Thus
\cite[Theorem~3.6]{matmat} gives continuous
orbit equivalence, hence the exact coefficient identities by the proof
of Theorem~\ref{thm:six-static}.

Each $A_j$ has zero diagonal, so $X_{E_j}$ has no fixed point.
The full $(j+1)$-shift has a fixed point. By the period inclusion in
\cite[Lemma~2.2]{shimomura}, conjugates of the full shift cannot
converge uniformly to the shift on $X_{E_j}$.
For $j=2$ the indegrees $(3,1,1)$
are all odd. Proposition~\ref{prop:six-ann} and iteration give the
inverse-image identity on every clopen set and at every iterate.
For $j=1$ the algebra is $\mathcal O_2$ by the diagonal isomorphism
above. Its gauge core has $K_0\cong\Z^3$: the stationary AF connecting
matrix $A_1^{\mathsf T}$ is unimodular ($\det A_1=1$), whereas the
binary full-shift core has $K_0=\Z[1/2]$ by the AF core description in \cite[Corollary~3.3 and the discussion following it]{raeburn} and \cite{cuntz1977}.
\end{exm}

\begin{exm}\label{ex:six-binary-trace}
A binary folding illustrates the extra constraint imposed by a trace
on the AF core. Let $G$ have states $1,2,3,4$, with letter $0$ sending $i$ to
$\max(1,i-1)$ and letter $1$ sending every state to $4$. Three letters
determine the endpoint: after the last $1$ only the trailing zeros
matter, and $000$ ends at $1$. The words $000,100,110,111$ reach states $1,2,3,4$, respectively. By Lemma~\ref{lem:folding-form},
$G$ is the quotient by a folding of $B(\{0,1\},3)$.
Proposition~\ref{prop:fold-matrix} identifies $C^*(G)$ with
$M_4(\mathcal O_2)$, which is isomorphic to $\mathcal O_2$ as a unital $C^*$-algebra
by \cite[Theorem~5.12]{aalp}.
The $K_0$- and $K_1$-groups of both algebras vanish.
Take a unital isomorphism and its inverse. The two factor maps of
Lemma~\ref{lem:fold-factors} are onto and every cylinder class is
zero. Corollary~\ref{cor:real-factor} applied to each constant
sequence therefore gives approximate $K$-conjugacy.

The two normalized core trace ranges are both $\Z[1/2]$ by
\eqref{eq:six-core}. In fact, a homeomorphism can match the two
traces on clopens: use the first two binary letters to choose the
initial state bijectively and read the remaining letters as outgoing
labels. A length-$n$ path cylinder then has preimage of binary measure
$2^{-(n+2)}$, exactly its core trace. Nevertheless no sequence can
match these clopen traces eventually and simultaneously approximate
the shifts from $\Sigma_2^+$ to the graph. Let $V_4$ be the cylinder
of paths starting at state $4$. Its core trace is $1/4$, while
$\sigma_G^{-1}V_4$ consists of four edge cylinders of total trace $4/8$.
Let $\nu$ be fair Bernoulli measure on $\Sigma_2^+$ and write
$\mu_G(D)=\tau_G(j_{X_G}(1_D))$ for clopen $D$.
Suppose $h_n:\Sigma_2^+\to X_G$ approximately intertwine the shifts
and $\nu(h_n^{-1}D)=\mu_G(D)$ eventually for each clopen $D$.
Choose $n$ large enough for this equality to hold for $V_4$ and
$\sigma_G^{-1}V_4$ and for
$h_n^{-1}(\sigma_G^{-1}V_4)=\sigma_2^{-1}(h_n^{-1}V_4)$,
as in Lemma~\ref{lem:real-transport}. Shift invariance of $\nu$ gives
\[
 \tfrac12=\mu_G(\sigma_G^{-1}V_4)
 =\nu(\sigma_2^{-1}(h_n^{-1}V_4))
 =\nu(h_n^{-1}V_4)=\mu_G(V_4)=\tfrac14,
\]
a contradiction.
\end{exm}

\subsection{A two-vertex shift beyond foldings}\label{sub:six-two}

We exhibit a shift outside the class of foldings that is approximately
$K$-conjugate to $\Sigma_3^+$. From the full shift to the edge shift
one factor map realizes the required classes at every length. In the
other direction we construct a sliding block code for each length $m$.

Consider the graph $E$ with adjacency matrix
$M=\left(\begin{smallmatrix}1&1\\2&0\end{smallmatrix}\right)$, in the
vertex order $A,B$. Denote its edges by
\[
 a:A\to A,\qquad b:A\to B,\qquad c_0,c_1:B\to A .
\]
The only loop is $a$, so $X_E$ has exactly one fixed point.
For any folding on $d\ge2$ letters the number of fixed points is
$\operatorname{tr}(M_G)=d$ by the characteristic polynomial in
Lemma~\ref{lem:six-det}.
Hence $X_E$ is not conjugate to the edge shift of any such folding.
The relations $[p_A]=[p_A]+[p_B]$ and $[p_B]=2[p_A]$
give $K_0(C^*(E))=\Z/2\Z$ with $[1]=[p_A]=\bar1$ and $[p_B]=\bar0$.
Also $\det(I-M)=-2$, so $K_1=0$. The graph is strongly connected, and both vertices emit two edges. As in Example~\ref{ex:six-loops},
\cite[Theorem~3.6]{matmat} gives a unital isomorphism with
$\mathcal O_3$. We identify both $K_0$-groups by their unit classes.
The indegrees $3,1$ are odd, so Proposition~\ref{prop:six-ann}
applies. We construct the required cylinder partitions in each direction.

\begin{prop}\label{prop:six-code}
There is a factor map $F:\Sigma_3^+\to X_E$ such that for every finite
path $\mu$,
\[
 c_{\Sigma_3^+}(F^{-1}([\mu]))=[p_{r(\mu)}]\quad\text{in }\Z/2\Z .
\]
\end{prop}
\begin{proof}
For $y\in\Sigma_3^+$ write $v_i=B$ when $y_i=2$ and
$y_{i+1}\in\{0,1\}$, and put $v_i=A$ otherwise. Consecutive vertices cannot both equal $B$. Define the $i$th edge of
$F(y)$ to be the unique edge from $v_i$ to $v_{i+1}$ if $v_i=A$,
and to be $c_{y_{i+1}}$ if $v_i=B$. In the latter case $v_{i+1}=A$,
so consecutive edges are composable. Each edge depends only on
$y_i,y_{i+1},y_{i+2}$. Hence $F$ is a sliding block code with
anticipation two.

To prove surjectivity, fix a path $x\in X_E$. Whenever $x_i=c_j$,
set $y_i=2$ and $y_{i+1}=j$. The pairs $\{i,i+1\}$ are disjoint
because no two consecutive edges have source $B$. At all other
positions set $y_i=0$. Then $v_i=B$ exactly when $x_i$ has source
$B$, and at each such position $F(y)_i=x_i$ by construction. At
positions with source $A$, the two endpoints determine the edge.
Thus $F(y)=x$.

For a finite path $\mu$ put $U_\mu=F^{-1}([\mu])$ and let
$x(\mu)=(x_0,x_1,x_2)$, where $x_i=c_{\Sigma_3^+}([i]\cap U_\mu)$
in $\Z/2\Z$. For the length-zero paths, $U_B=[20]\sqcup[21]$
and $U_A=\Sigma_3^+\setminus U_B$. Each word cylinder has class
$\bar1$ by Corollary~\ref{cor:fold-scaled-k}, so
$x(A)=(1,1,1)$ and $x(B)=(0,0,0)$.
Put $u_A=(1,1,1)$, $u_B=(0,0,1)$ and
$\epsilon_A=1$, $\epsilon_B=0$. We claim that
\begin{equation}\label{eq:six-parity}
 x(\mu)=\epsilon_{r(\mu)} u_{s(\mu)}.
\end{equation}
Write $iD=\{iy:y\in D\}$ for a letter $i$ and a clopen set $D$.
If $r(e)=s(\mu)$, the definition of $F$ gives
\[
 U_{e\mu}=
 \begin{cases}
  0U_\mu\sqcup 1U_\mu\sqcup 2([2]\cap U_\mu), &e=a,b,\\
  2([j]\cap U_\mu), &e=c_j.
 \end{cases}
\]
Indeed the initial vertex is $A$ precisely on $[0]\sqcup[1]\sqcup[22]$,
and the initial edge is $c_j$ precisely on $[2j]$. The initial vertex
of $F(y)$ specifies the range of the edge prepended to $F(y)$.
Since $c_{\Sigma_3^+}(iD)=c_{\Sigma_3^+}(D)$ by the cylinder formula
and finite additivity, putting $S=x_0+x_1+x_2$ gives
\[
 x(e\mu)=
 \begin{cases}
 (S,S,x_2),&e=a,b,\\
 (0,0,x_j),&e=c_j.
 \end{cases}
\]
Both $u_A$ and $u_B$ have
coordinate sum $1$ and last coordinate $1$. Hence the first rule
sends either row to $u_A$, and the second sends $u_A$ to $u_B$.
Induction gives \eqref{eq:six-parity}. Summing its coordinates proves
the assertion, including paths ending at $B$ whose nonempty preimages
have zero class.
\end{proof}

\begin{prop}\label{prop:six-marker}
For each $m\ge2$ there is a continuous shift-commuting map
$G_m:X_E\to\Sigma_3^+$ such that
$c_{X_E}(G_m^{-1}([w]))=\bar1$ for every ternary word $w$ of length
$m$. In particular every one of these cylinders meets the image.
\end{prop}
\begin{proof}
Use the Euler circuit construction of a de Bruijn sequence
\cite[Chapter~8]{vlw} on the ternary de Bruijn graph of order $m-1$,
starting at the vertex $0^{m-1}$. Reading the initial vertex and the
successive edge labels gives a word $W=W_0\cdots W_{T-1}$ with
$T=3^m+m-1$. Every ternary word of length $m$ occurs exactly once
in $W$, and the prefix and suffix of length $m-1$ both equal $0^{m-1}$.

Take $N=T+2m$ and the path cylinder $U=[c_0a^N]\subset X_E$.
Two occurrences of this path have starting times more than $N$ apart:
the $N$ edges after $c_0$ are all $a$. The path ends at $A$, so
$c_{X_E}(U)=\bar1$. Proposition~\ref{prop:six-ann} gives
$c_{X_E}(\sigma_E^{-j}U)=\bar1$ for every $j\ge0$.
Define a continuous letter function $g:X_E\to\{0,1,2\}$ by
\[
 g(x)=W_{T-1-j}\ \text{if }\sigma_E^j x\in U\text{ for }0\le j<T,
 \qquad g(x)=0\ \text{otherwise}.
\]
The indicated clopen sets are disjoint. Set
$G_m(x)_i=g(\sigma_E^i x)$ for $i\ge0$.

For $w\ne0^m$, let $k$ be its unique starting position in $W$.
We have
\begin{equation}\label{eq:six-marker}
 G_m^{-1}([w])=\sigma_E^{-(T-1-k)}U .
\end{equation}
Call $t\ge0$ a marker time for $x$ if $\sigma_E^t x\in U$.
For every such $t$,
$G_m(x)_i=W_{i-(t-T+1)}$ whenever $i\in[t-T+1,t]\cap\mathbb N_0$.
These intervals are disjoint by the separation of marker times, and
$G_m(x)_i=0$ outside their union.

For the inclusion from left to right, take $G_m(x)\in[w]$ and choose
$i\in\{0,\ldots,m-1\}$ with $w_i\ne0$. There is a unique marker time $t$ with
$i\in[t-T+1,t]$. The index $j=i-(t-T+1)$ satisfies
$m-1\le j\le T-m$ by the two zero blocks at the ends of $W$. Thus
$t-T+1=i-j\le0$ and $t=i-j+T-1\ge m-1$. It follows that $[0,m-1]\subseteq[t-T+1,t]$.
Thus $w$ begins at position $T-1-t$ in $W$. Its unique occurrence
begins at $k$, so $t=T-1-k$.
Conversely, if $t=T-1-k$ is a marker time, then for $0\le i<m$
we have $0\le t-i<T$ since $0\le k\le T-m$. The definition of
$g$ gives $G_m(x)_i=W_{T-1-(t-i)}=W_{k+i}=w_i$.
This proves \eqref{eq:six-marker}.
Thus the preimages of all $3^m-1$ nonzero words have class $\bar1$.
The remaining preimage has class
$[1]-(3^m-1)\bar1=\bar1$ in $\Z/2\Z$. Every preimage is therefore
nonempty.
\end{proof}

\begin{cor}\label{cor:six-two}
$X_E$ and $\Sigma_3^+$ are approximately $K$-conjugate.
\end{cor}
\begin{proof}
Fix a unital $*$-isomorphism $\Psi:\mathcal O_3\to C^*(E)$.
It induces the identification by unit classes on $K_0$. Both
shift spaces are Cantor spaces, since every nonempty cylinder admits
at least two distinct one-edge extensions. Both algebras are unital,
simple and purely infinite. Apply Corollary~\ref{cor:real-factor}
in the following two directions.
\begin{enumerate}
\item For $Z=X_E$, $W=\Sigma_3^+$, use $F_n=G_{n+1}$ and $\theta=\Psi$.
Proposition~\ref{prop:six-marker} gives nonempty inverse images for
every length-$(n+1)$ cylinder, each of class $\bar1=\Psi_{*0}(\bar1)$
under our identification.
\item For $Z=\Sigma_3^+$, $W=X_E$, use the constant sequence $F_n=F$
and $\theta=\Psi^{-1}$. Surjectivity gives nonempty inverse images,
including for cylinders ending at $B$. Proposition~\ref{prop:six-code}
gives their classes, including the zero class for those ending at $B$.
\end{enumerate}
The two applications give the two uniform limits in \eqref{eq:ak-uniform}
and the two norm limits in \eqref{eq:ak-norm} along the same respective
homeomorphisms. Their algebra isomorphisms have fixed inverse
$KL$-classes $[\Psi]_{KL}$ and $[\Psi^{-1}]_{KL}$.
\end{proof}

The two shifts in Corollary~\ref{cor:six-two} have different topological
entropies. There are $2\cdot 2^n$ edge paths of length $n$ in $E$ and
$3^n$ ternary words. The word-growth formula for entropy
\cite[Chapter~4]{lindmarcus} gives $\log 2$ and $\log 3$, respectively.
In particular, $X_E$ and $\Sigma_3^+$ are not conjugate, and no
single factor map from $X_E$ onto $\Sigma_3^+$ can replace the
sequence $(G_m)$, since entropy cannot increase under a factor map
\cite[Chapter~4]{lindmarcus}.
For a general primitive graph it remains to determine whether a
unital algebra isomorphism with $\Od$, the indegree congruences and
the existence of a fixed point suffice for approximate $K$-conjugacy.
The constructions of $F$ and $G_m$ use the specific transitions of $E$.


\begin{thebibliography}{99}
\bibitem{aalp} G. Abrams, P.~N. {\'A}nh, A. Louly and E. Pardo, The classification question for Leavitt path algebras, \textit{J. Algebra} \textbf{320} (2008), 1983--2026.
\bibitem{abramssmith} G. Abrams and C. Smith, The matrix type of purely infinite simple Leavitt path algebras, \textit{Acta Math. Hungar.} \textbf{128} (2010), 381--385.
\bibitem{blackadar} B. Blackadar, \textit{$K$-Theory for Operator Algebras}, 2nd ed., Mathematical Sciences Research Institute Publications, vol.~5, Cambridge University Press, Cambridge, 1998.
\bibitem{bleak} C. Bleak, P.~J. Cameron and F. Olukoya, Automorphisms of shift spaces and the Higman--Thompson groups: the one-sided case, \textit{Discrete Anal.} (2021), Paper No.~15, 35 pp.
\bibitem{brixconj} K.~A. Brix and T.~M. Carlsen, Cuntz--Krieger algebras and one-sided conjugacy of shifts of finite type and their groupoids, \textit{J. Aust. Math. Soc.} \textbf{109} (2020), 289--298.
\bibitem{brixcovers} K.~A. Brix and T.~M. Carlsen, $C^*$-algebras, groupoids and covers of shift spaces, \textit{Trans. Amer. Math. Soc. Ser. B} \textbf{7} (2020), 134--185.
\bibitem{bcw} N. Brownlowe, T.~M. Carlsen and M.~F. Whittaker, Graph algebras and orbit equivalence, \textit{Ergodic Theory Dynam. Systems} \textbf{37} (2017), 389--417.
\bibitem{carlsencp} T.~M. Carlsen, Cuntz--Pimsner $C^*$-algebras associated with subshifts, \textit{Internat. J. Math.} \textbf{19} (2008), 47--70.
\bibitem{carlsenrout} T.~M. Carlsen and J. Rout, Diagonal-preserving gauge-invariant isomorphisms of graph $C^*$-algebras, \textit{J. Funct. Anal.} \textbf{273} (2017), no.~9, 2981--2993.
\bibitem{carlsensilvestrov} T.~M. Carlsen and S. Silvestrov, $C^*$-crossed products and shift spaces, \textit{Expo. Math.} \textbf{25} (2007), no.~4, 275--307. 
\bibitem{cuntz1977} J. Cuntz, Simple $C^*$-algebras generated by isometries, \textit{Comm. Math. Phys.} \textbf{57} (1977), 173--185.
\bibitem{cuntz} J. Cuntz, $K$-theory for certain $C^*$-algebras, \textit{Ann. of Math.} (2) \textbf{113} (1981), 181--197.
\bibitem{ck} J. Cuntz and W. Krieger, A class of $C^*$-algebras and topological Markov chains, \textit{Invent. Math.} \textbf{56} (1980), 251--268.
\bibitem{dadarlat} M. Dadarlat, Approximately unitarily equivalent morphisms and inductive limit $C^*$-algebras, \textit{$K$-Theory} \textbf{9} (1995), 117--137.
\bibitem{gps} T. Giordano, I.~F. Putnam and C.~F. Skau, Topological orbit equivalence and $C^*$-crossed products, \textit{J. Reine Angew. Math.} \textbf{469} (1995), 51--111.
\bibitem{higman} G. Higman, \textit{Finitely Presented Infinite Simple Groups}, Notes on Pure Mathematics, no.~8, Australian National University, Canberra, 1974.
\bibitem{kechris} A.~S. Kechris, \textit{Classical Descriptive Set Theory}, Graduate Texts in Mathematics, vol.~156, Springer-Verlag, New York, 1995.
\bibitem{llw} H. Lin, H. Liu and S. Wei, Approximate $K$-conjugacy of minimal Cantor--circle extensions I: The rigid case, manuscript.
\bibitem{linmatui} H. Lin and H. Matui, Minimal dynamical systems and approximate conjugacy, \textit{Math. Ann.} \textbf{332} (2005), 795--822.
\bibitem{lindmarcus} D. Lind and B. Marcus, \textit{An Introduction to Symbolic Dynamics and Coding}, 2nd ed., Cambridge University Press, Cambridge, 2021.
\bibitem{matuifull} H. Matui, Topological full groups of one-sided shifts of finite type, \textit{J. Reine Angew. Math.} \textbf{705} (2015), 35--84.
\bibitem{matmat} K. Matsumoto and H. Matui, Continuous orbit equivalence of topological Markov shifts and Cuntz--Krieger algebras, \textit{Kyoto J. Math.} \textbf{54} (2014), 863--877.
\bibitem{phillips} N.~C. Phillips, A classification theorem for nuclear purely infinite simple $C^*$-algebras, \textit{Doc. Math.} \textbf{5} (2000), 49--114.
\bibitem{putnam} I.~F. Putnam, The $C^*$-algebras associated with minimal homeomorphisms of the Cantor set, \textit{Pacific J. Math.} \textbf{136} (1989), 329--353.
\bibitem{raeburn} I. Raeburn, \textit{Graph Algebras}, CBMS Regional Conference Series in Mathematics, vol.~103, American Mathematical Society, Providence, RI, 2005.
\bibitem{rll} M. R{\o}rdam, F. Larsen and N.~J. Laustsen, \textit{An Introduction to $K$-Theory for $C^*$-Algebras}, London Mathematical Society Student Texts, vol.~49, Cambridge University Press, Cambridge, 2000.
\bibitem{rs} J. Rosenberg and C. Schochet, The K{\"u}nneth theorem and the universal coefficient theorem for Kasparov's generalized $K$-functor, \textit{Duke Math. J.} \textbf{55} (1987), 431--474.
\bibitem{shimomura} T. Shimomura, Chain mixing endomorphisms are approximated by subshifts on the Cantor set, \textit{Tsukuba J. Math.} \textbf{35} (2011), no.~1, 67--77.
\bibitem{vlw} J.~H. van Lint and R.~M. Wilson, \textit{A Course in Combinatorics}, 2nd ed., Cambridge University Press, Cambridge, 2001.
\bibitem{webster} S.~B.~G. Webster, The path space of a directed graph, \textit{Proc. Amer. Math. Soc.} \textbf{142} (2014), 213--225.
\bibitem{weihe} S. Wei and Z. He, Dimensions associated with surjective local homeomorphisms and subshifts with low complexity, \textit{Proc. Roy. Soc. Edinburgh Sect. A} (2026), 1--34.
\end{thebibliography}
\end{document}